\documentclass[11pt]{amsart}

\usepackage{amsmath}
\usepackage{amscd}
\usepackage{amsthm}
\usepackage{hyperref}
\usepackage[margin=1.0in]{geometry}

\numberwithin{equation}{section}

\newtheorem{prop}{Proposition}
\newtheorem{lemma}[prop]{Lemma}
\newtheorem{thm}[prop]{Theorem}
\newtheorem{cor}[prop]{Corollary}
\newtheorem{conj}[prop]{Conjecture}
\numberwithin{prop}{section}

\theoremstyle{definition}
\newtheorem{defn}[prop]{Definition}

\newcommand{\dt}{\frac{\partial}{\partial t}}

\newcommand{\brs}[1]{\left| #1 \right|}

\renewcommand{\gg}{\gamma}
\newcommand{\gD}{\Delta}
\newcommand{\gd}{\delta}
\newcommand{\gs}{\sigma}

\newcommand{\gl}{\lambda}
\newcommand{\gt}{\theta}
\newcommand{\gw}{\omega}
\newcommand{\ga}{\alpha}
\newcommand{\gb}{\beta}
\renewcommand{\ge}{\epsilon}
\newcommand{\N}{\nabla}

\renewcommand{\bar}[1]{\overline{#1}}
\newcommand{\del}{\partial}
\newcommand{\delb}{\bar{\partial}}
\newcommand{\hi}{\frac{\sqrt{-1}}{2}}

\newcommand{\bi}{\bar{i}}
\newcommand{\bj}{\bar{j}}
\newcommand{\bk}{\bar{k}}
\newcommand{\bl}{\bar{l}}
\newcommand{\bm}{\bar{m}}
\newcommand{\bn}{\bar{n}}

\newcommand{\bq}{\bar{q}}
\newcommand{\br}{\bar{r}}
\newcommand{\bs}{\bar{s}}

\newcommand{\bv}{\bar{v}}

\newcommand{\U}{\Upsilon}

\newcommand{\ra}{\rightarrow}
\renewcommand{\to}{\tilde{\omega}}

\newcommand{\til}[1]{\widetilde{#1}}

\newcommand{\seven}{\mbox{VII}}

\newcommand{\LL}{\mathcal L}

\newcommand{\hook}{\mathbin{\hbox{\vrule height2.4pt width4.5pt depth-2pt
\vrule height5pt width0.4pt depth-2pt}}}

\DeclareMathOperator{\Rc}{Rc}
\DeclareMathOperator{\Ric}{Ric}

\DeclareMathOperator{\inj}{inj}
\DeclareMathOperator{\tr}{tr}
\DeclareMathOperator{\Ker}{Ker}
\DeclareMathOperator{\Img}{Im}

\DeclareMathOperator{\Vol}{Vol}

\DeclareMathOperator{\End}{End}
\DeclareMathOperator{\Diff}{Diff}
\DeclareMathOperator{\Met}{Met}

\begin{document}

\title[Regularity of Pluriclosed Flow]{Regularity Results for Pluriclosed
Flow}

\begin{abstract} In \cite{ST2} the authors introduced a parabolic flow of
pluriclosed metrics.  Here we give improved regularity results for solutions to
this equation.  Furthermore, we exhibit this equation as the gradient flow of
the lowest eigenvalue of a certain Schr\"odinger operator, and show the
existence of an expanding entropy functional for this flow.  Finally, we
motivate a conjectural picture of the optimal
regularity results for this flow, and discuss some of the consequences.
\end{abstract}

\author{Jeffrey Streets}
\thanks{Streets supported by DMS-1201569, Tian supported by DMS-0804095}
\address{Rowland Hall\\
         University of California, Irvine\\
         Irvine, CA 92617}
\email{\href{mailto:jstreets@uci.edu}{jstreets@uci.edu}}

\author{Gang Tian}
\address{Fine Hall\\
	 Princeton University\\
	 Princeton, NJ 08544}
\email{\href{mailto:tian@math.princeton.edu}{tian@math.princeton.edu}}

\date{February 19th 2013}


\maketitle

\section{Introduction}

Let $(M^{2n}, J)$ be a complex manifold, and let $\omega$ denote a Hermitian
metric on $M$.  The metric $\omega$ is \emph{pluriclosed} if
\begin{align*}
 \del \delb \omega = 0.
\end{align*}
\noindent Consider the initial value problem
\begin{gather} \label{PCF}
\begin{split} 
\dt \omega =&\ \del \del^* \omega + \delb \delb^* \omega + \frac{\sqrt{-1}}{2}
\del \delb \log \det g\\
\omega(0) =&\ \omega_0.
\end{split}
\end{gather}
\noindent This equation was introduced in \cite{ST2} as a tool for understanding
complex,
non-K\"ahler manifolds.  Equation (\ref{PCF}) falls into a general class of
flows of Hermitian
metrics, and as shown in \cite{ST1}, solutions to (\ref{PCF}) exist as long as
the Chern curvature, torsion, and covariant derivative of torsion are bounded.
This is analogous to the long time existence theorem by R. Hamilton (\cite{Ham3}
Theorem 14.1) which states that
the Ricci flow with any initial data has a solution on $[0,T)$, where either $T
= \infty$ or the curvature of the solution blows up at time $T$. A natural
problem is whether or not we can drop the hypothesis that the torsion and its
first covariant derivative is bounded at a finite singular time for (\ref{PCF}).
 In the case $n = 2$, we showed in \cite{ST2} that a bound on the
Chern curvature suffices to show long time existence. The difficulty in general
arises from the fact that the induced evolution equation 
on the Chern curvature involves the torsion and its derivatives.  Our crucial
observation for overcoming this difficulty is that the Bismut 
connection is a much more natural connection for studying (\ref{PCF}).

In this paper, we will give sharper long time existence theorems for
(\ref{PCF}).  We will also prove some useful regularity theorems and paint a
more concrete picture for the conjectural existence and singularity formation
for the flow (\ref{PCF}).  Furthermore, by giving an interesting interpretation
of the
flow using the Bismut connection, we exhibit a remarkable relationship of
(\ref{PCF}) to mathematical physics.  Specifically, we show that up to gauge
equivalence (\ref{PCF}) is the renormalization group flow of a nonlinear sigma
model with nonzero $B$-field.  As a consequence we derive that (\ref{PCF}) is a
gradient flow, and exhibit a certain entropy functional.  Finally, we discuss
some applications of our conjectural picture to understanding the topology of
Class $\seven$ surfaces.

We start by recalling the \emph{Bismut connection}.  Let $(M^{2n}, \omega, J)$
be a complex
manifold with pluriclosed metric.  Let $D$ denote the Levi Civita connection. 
Then the Bismut connection $\N$ is defined via
\begin{align*}
\left< \N_X Y, Z \right> = \left<D_X Y, Z \right> + \frac{1}{2} d^c \omega(X, Y,
Z)
\end{align*}
where $d^c \omega(X, Y, Z) := d \omega(JX, JY, JZ)$.  Let $\Omega$ denote the
curvature of this connection, and let $P$ denote the
Chern form of this connection, i.e. in complex coordinates
\begin{align*}
P_{i \bj} =&\ \Omega_{i \bj k}^k
\end{align*}
Finally, let $P^C$ denote the Ricci form associated to the Chern connection. 
One can calculate (\cite{AlexIvan}) that
\begin{align*}
P =&\ P^C - d d^* \omega.
\end{align*}
In particular, this implies that a solution to (\ref{PCF}) may be expressed as
\begin{align} \label{PCFBismut}
\dt \omega =&\ - P^{1,1}
\end{align}
where $P^{1,1}$ denotes the projection of $P$ onto $(1,1)$-forms.  This is a
convenient framework for understanding solutions to (\ref{PCF}).  In particular,
with the clarifying lens of this connection, we are able to show that
(\ref{PCF}) is the gradient flow of the first eigenvalue of a particular
Schr\"odinger operator.  First generalize the notation slightly and let $(M^n,
g)$ be a Riemannian manifold, and let $T$ denote
a three-form on $M$.  Let
\begin{align*}
\mathcal F(g, T, f) =&\ \int_M \left[ R - \frac{1}{12} \brs{T}^2 + \brs{\N f}^2
\right] e^{-f} dV.
\end{align*}
Furthermore set
\begin{align*}
\gl(g, T) = \inf_{\{f | \int_M e^{-f} dV = 1 \}} \mathcal F(g, T, f).
\end{align*}

\noindent In section \ref{gradient} we exhibit equation (\ref{PCF}) as the
gradient flow of $\gl$.

\begin{thm} \label{gradientPCF} Let $(M^{2n}, \omega, J)$ be a complex manifold
with pluriclosed metric.  Let $\omega(t)$ denote the solution to (\ref{PCF})
with initial condition $\omega$, and let $g(t)$ be the associated metric, and
$T(t)$ the torsions of the associated Bismut connections.  Let $\Met$ denote the
space of smooth metrics on $M$, and let
\begin{align*}
\mathcal M := \frac{\{(g, T) | g \in \Met, T \in \Lambda^3, dT = 0
\}}{\Diff_+(M)}
\end{align*}
where $\Diff_+$ is the group of oriented diffeomorphisms of $M$, acting
naturally on $g$ and $T$.  The pair $(g(t), T(t))$ is a solution of the gradient
flow of $\gl$ acting on $\mathcal M$.
\end{thm}

More specifically, we show that after pulling back by the one-parameter family
of
diffeomorphisms generated by the Lee forms of the time dependent metrics,
equation (\ref{PCF}) is unmasked as the $B$-field renormalization group flow of
string theory.  This flow has been previously studied, and admits the
generalization $\gl$ of the Perelman energy (\cite{Woolgar}).  Furthermore an
expanding entropy functional for this flow was discovered by the first named
author (\cite{Streets}), which hence is monotone for solutions to (\ref{PCF}) as
well.  These observations show that any breather solution is automatically a
gradient soliton (Corollary \ref{nobreathers}), and furthermore imply strong
results on certain long-time solutions.

Turning to the regularity theory, we show that a bound on the Bismut Ricci
curvature suffices to obtain long time existence for solutions of (\ref{PCF}).

\begin{thm} \label{ricciregularity} Let $(M^{2n}, \omega(t), J)$ be a solution
to (\ref{PCF}) on $[0, \tau)$.  Suppose
\begin{align*}
\int_0^{\tau} \sup_{M \times \{t\}} \brs{P^{1,1}} dt < \infty.
\end{align*}
Then the solution extends smoothly past time $\tau$.
\end{thm}
\noindent This theorem is analogous to a result of N. Sesum for the Ricci flow
(\cite{Sesum} Theorem 2), and already represents a significant improvement, as
we have reduced the
regularity requirement to understanding the Ricci-type curvature of a specific
connection.

The theory of K\"ahler Ricci flow is considerably more developed than the
general study of Ricci flow.  One of key reasons for this is the reduction of
the K\"ahler-Ricci flow to a scalar equation.  Inspired by this, we will
introduce a certain potential function $\phi$ along a solution to (\ref{PCF})
and prove a regularity theorem in term of this potential and the torsion.
Let $(M^{2n}, \til{\omega}, J)$ be a complex manifold with pluriclosed metric,
and let $\omega(t)$ be a solution to (\ref{PCF}).  We define

\begin{gather} \label{potentialeqn}
\begin{split}
\dt \phi - \gD \phi =&\ \tr_{\gw} \to - n\\
\phi(0) =&\ 0.
\end{split}
\end{gather}
Here $\gD$ is the canonical Laplacian associated to the time dependent metric
$\omega(t)$, i.e. $\gD = \tr_{\omega} \del \delb$.  It follows from standard
parabolic theory that $\phi$ exists on the same time interval that $\omega(t)$
exists.  More generally, one may define this with respect to a one-parameter
family of background metrics $\til{\omega}(t)$.  Related quantities for
Monge-Ampere type equations on almost K\"ahler surfaces were considered in
\cite{Whinekove}.

\begin{thm} \label{regularity2} Let $(M^{2n}, \til{g}, J)$ be a compact complex
manifold and suppose $g(t)$ is a solution to (\ref{PCF})
on $[0, \tau)$ and suppose there is a constant $C$ such that
\begin{align*}
\sup_{M \times [0, \tau)} \brs{\phi} \leq&\ C,\\
\sup_{M \times [0, \tau)} \brs{T}^2 \leq&\ C.
\end{align*}
Then $g(t) \rightarrow g(\tau)$ in $C^{\infty}$, and the flow extends smoothly
past time $\tau$.
\end{thm}

\noindent This represents a significant reduction of the regularity requirement
for solutions to (\ref{PCF}), effectively reducing the question to understanding
the behavior of the potential function and the torsion.  The proof
involves
applying the maximum principle to carefully chosen quantities.  Going further,
one would like to understand what the optimal existence and
regularity theorems are for (\ref{PCF}).  For this we again take a cue from the
study
of K\"ahler Ricci flow.  Suppose $(M^{2n}, \omega_0, J)$ is a K\"ahler manifold,
and recall the K\"ahler Ricci flow equation
\begin{gather} \label{KRF}
\begin{split}
 \dt \omega =&\ \frac{\sqrt{-1}}{2} \del \delb \log \det g\\
\omega(0) =&\ 0.
\end{split}
\end{gather}
Associated to a solution $\omega(t)$ of (\ref{KRF}) is an ODE in $H^2(M, \mathbb
R)$ which has solution
\begin{align*}
 [\omega(t)] = [\omega_0] - t c_1(M).
\end{align*}
The optimal regularity theorems for K\"ahler Ricci flow assert that as long as
the solution to the ODE above remains in the K\"ahler cone, the solution exists
up to that time (\cite{TZ}, see also \cite{Tian07}).  An essential ingredient of
these
theorems is the reduction of K\"ahler Ricci flow to a scalar equation, then
exploiting the estimates of the Monge Ampere equation using the maximum
principle.

Keeping with this theme, observe that a pluriclosed metric defines a class in
the Aeppli cohomology group
\begin{align*}
\mathcal H^{1,1}_{\del + \delb} = \frac{ \{ \Ker \del \delb : \Lambda_{\mathbb
R}^{1,1} \rightarrow \Lambda_{\mathbb R}^{2,2} \} }{ \{ \del \ga + \delb
\bar{\ga} | \ga \in \Lambda^{0,1} \} }.
\end{align*}
Define the space
$\mathcal P_{\del + \delb}$ to be the cone of the classes in $\mathcal
H_{\del + \delb}$ which contain positive definite elements.  Solutions to
(\ref{PCF}) clearly define ODE's in $\mathcal H^{1,1}_{\del + \delb}$, and it is
natural to conjecture (see Conjecture \ref{existenceconj}) that the maximal
existence time is characterized by the first time at which the boundary of
$\mathcal P_{\del + \delb}$ is reached.  If true, this would have strong
implications on complex surfaces since the cone $\mathcal P_{\del + \delb}$ is
essentially characterized on complex surfaces in terms of the action of the
class in $\mathcal H_{\del + \delb}$ on curves.  More precisely, {\it if $\phi$
is a pluriclosed
(1,1)-form on a complex non-K\"ahler surface $(M^4, J)$, then $\phi \in \mathcal
P_{\del + \delb}$ if and
only if (1) $\int_M \phi \wedge \gamma_0 > 0$; (2) $\int_D \phi > 0$ for every
effective divisor with negative self
intersection.}  Here $\gamma_0$ is the kernel of the projection map from the
$(1,1)$ Bott-Chern cohomology of $M$ to $H^{1,1}$, explained further in section
\ref{conesection}.

To further illustrate the significance of the cone $\mathcal P_{\del + \delb}$,
we show
in section \ref{conesection} that as long as the solution to the associated ODE
remains in the
interior of $\mathcal P$, solutions to (\ref{PCF}) on complex surfaces may be
canonically reduced to solutions of a certain PDE on $\ga \in
\Lambda^{0,1}$, and an ODE on $\psi \in \Lambda^{1,1}$.  Specifically, we can
find a background metric $\til{g}$ so
that, setting $\omega(t) = \omega(0) + \del \alpha + \delb \bar{\ga} + \psi$,
one has the solution to (\ref{PCF}) reduced to
\begin{gather} \label{aeqn}
\begin{split}
\dt \ga =&\ \del^{*}_{\gw} \gw + \frac{\sqrt{-1}}{4} \delb \log
\frac{\gw^{n}}{\to^{n}},\\
\dt \psi =&\ - c_1(\til{g})\\
\ga(0) =&\ \ga_0, \qquad \psi(0) = \psi_0.
\end{split}
\end{gather}
This equation may be taken as an ansatz for equation (\ref{PCF}) in \emph{any}
dimension.  We discuss some further properties of this equation in section
\ref{conesection}.

Finally, in section \ref{nonsingularsection} we examine nonsingular solutions to
(\ref{PCF}), suitably normalized, on Class $\seven^+$ surfaces.  By exploiting
certain monotone
quantities we show that a positive resolution of the
conjectural picture of singularity
formation outlined in section \ref{conesection} implies the existence of a curve
on a Class $\seven^+$ surface. 
Due to the results of Nakamura \cite{Nakamura}, and Dloussky, Oeljeklaus, and
Toma
\cite{DOT}, the classification of such surfaces reduces to finding sufficiently
many curves.  In particular we show that our conjectural picture implies the
classification of Class $\seven^+$ surfaces with $b_2 = 1$.  This is discussed
further in section \ref{nonsingularsection}.

Here is an outline of the rest of the paper.  In section 2 we establish certain
differential inequalities for solutions to (\ref{PCF}) which are inspired by the
theory of the complex Monge-Amp\`ere equation which can be used to establish
uniform bounds on the metric.  Next in section 3 we derive a $C^1$ estimate for
the metric along solutions to (\ref{PCF}) under certain hypotheses.  Building on
these estimates, in section 4 we give the proofs of Theorems
\ref{ricciregularity} and \ref{regularity2}.  In section 5 we
outline a conjectural picture of formation of singularities of solutions to
(\ref{PCF}) in any dimension, and further clarify this picture in the case of
complex surfaces.  In section \ref{gradient} we show that (\ref{PCF}) is a
gradient flow, and show the existence of an entropy functional.  In section
\ref{nonsingularsection} we derive some
consequences of the conjectural regularity picture developed in section
\ref{conesection}, and
section \ref{conclusion} is a brief
conclusion.

The first named author would like to thank his advisor, Mark Stern, for first
introducing him to the $B$-field renormalization group flow, and for many
interesting conversations.  He would also like to thank Nicholas Buchdahl for
interesting conversations.

\section{\texorpdfstring{$L^{\infty}$}{Uniform} metric estimate}

In this section we derive differential inequalities and produce a priori
estimates for the metric along a
solution to (\ref{PCF}) which are similar in spirit to the Laplacian estimates
for the solution to the complex Monge Amp\`ere equation.  We will fix a
one-parameter family of background metrics $\til{\omega}(t)$ and assume that
they are uniformly bounded on the time interval of consideration.  Furthermore,
we set
\begin{align*}
\psi := \dt \til{\omega}.
\end{align*}
Finally, $\phi$ will always denote the solution to (\ref{potentialeqn}) taken
with respect to $\til{\omega}(t)$.

\begin{lemma} \label{Phicoordinates} Let $(M^{2n}, \omega, J)$ be a complex
manifold with pluriclosed metric.  Then in local complex coordinates,
\begin{align*}
P^{1,1}(\omega)_{k \bl} =&\ g^{i \bj} \left( - g_{k \bl, i \bj} + g^{m \bn}
g_{k \bn, i} g_{\bl m, \bj} \right) - g^{m \bn} g^{p \bq} \left( g_{m \bq, k} -
g_{k \bq, m} \right) \left( g_{\bn p, \bl} - g_{\bl p, \bn} \right).
\end{align*}
\begin{proof} As exhibited in \cite{ST2}, one has an expression
\begin{align*}
P^{1,1}(\omega)_{k \bl} =&\ - S_{k \bl} + Q^1_{k \bl}
\end{align*}
where
\begin{align*}
S_{k \bl} =&\ g^{i \bj} \Omega_{i \bj k \bl}
\end{align*}
and
\begin{align*}
Q^1_{k \bl} =&\ g^{m \bn} g^{p \bq} T_{k m \bq} T_{\bl \bn p},
\end{align*}
where $\Omega$ and $T$ are the curvature and torsion of the Chern connection,
respectively.  The lemma then follows from direct calculations.
\end{proof}
\end{lemma}

\begin{lemma} \label{prelimcalcs} Let $(M^{2n}, \til{\omega}(t), J)$ be a
complex
manifold with a one-parameter family of pluriclosed metrics.  Let $\gw(t)$
denote a solution to
(\ref{PCF}).  Then in local complex coordinates,
\begin{align*}
\dt \tr_{\to} \gw =&\ \til{g}^{k \bl} g^{i \bj} \left( g_{k\bl, i \bj}
\right) - \til{g}^{k \bl} g^{i \bj} g^{m \bn} g_{k \bn, i}
g_{\bl m, \bj}\\
&\ + \til{g}^{k \bl} g^{m \bn} g^{p \bq} \left( g_{m \bq, k} -
g_{k \bq, m} \right) \left( g_{\bn p, \bl} - g_{\bl p, \bn}
\right) - \left<\psi, \omega \right>_{\to},\\
\dt \log \frac{\gw^n}{\to^n} =&\ \gD \log \frac{\gw^n}{\to^n} +
g^{k \bl} g^{m \bn} g^{p \bq} \left( g_{m \bq, k} -
g_{k \bq, m} \right) \left( g_{\bn p, \bl} - g_{\bl p, \bn}
\right)\\
&\ + \tr_{\gw} \del \delb \log \det \til{g} - \tr_{\to} \psi,\\
\dt \tr_{\gw} \to =&\ g^{i \bq} g^{p \bj} g^{r \bs} \left[ - g_{p \bq, r \bs} +
g^{u \bv} g_{p \bv, r} g_{u \bq, \bs} \right] \til{g}_{i \bj}\\
&\ - g^{i \bq} g^{p \bj} g^{r \bs} g^{u \bv} \left( g_{r \bv, p} - g_{p \bv, r}
\right) \left( \del_{\bq} g_{u \bs} - \del_{\bs} g_{u \bq} \right) \til{g}_{i
\bj} + \tr_{\gw} \psi.
\end{align*}
\begin{proof} Starting from Lemma \ref{Phicoordinates}, we compute
\begin{align*}
\dt \tr_{\to} \gw =&\ - \til{g}^{k \bl} P^{1,1} (\gw)_{k \bl} - \left<\psi,
\omega
\right>_{\to}\\
=&\ \til{g}^{k \bl} g^{i \bj} \left( g_{k\bl, i \bj} \right) - \til{g}^{k \bl}
g^{i \bj} g^{m \bn} g_{k \bn, i} g_{\bl m, \bj}\\
&\ + \til{g}^{k \bl} g^{m \bn} g^{p \bq} \left( g_{m \bq, k} -
g_{k \bq, m} \right) \left( g_{\bn p, \bl} - g_{\bl p, \bn}
\right) - \left<\psi, \omega \right>_{\to}.
\end{align*}
Using a general calculation we have $\dt \log \frac{\gw^n}{\to^n} = \tr_{\gw}
\left(\dt \gw \right) - \tr_{\to} \left( \dt \to \right)$.  Thus
\begin{align*}
\dt \log \frac{\gw^n}{\to^n} =&\ g^{k \bl} g^{i \bj} g_{k \bl,
i \bj} - g^{k \bl} g^{i \bj} g^{m \bn} g_{k \bn, i}
g_{\bl m, \bj}\\
&\ + g^{k \bl} g^{m \bn} g^{p \bq} \left( g_{m \bq, k} -
g_{k \bq, m} \right) \left( g_{\bn p, \bl} - g_{\bl p, \bn}
\right) - \tr_{\to} \psi\\
=&\ \gD \log \gw^n + g^{k \bl} g^{m \bn} g^{p \bq}
\left( g_{m \bq, k} - g_{k \bq, m} \right) \left( g_{\bn p,
\bl} - g_{\bl p, \bn} \right) - \tr_{\to} \psi\\
=&\ \gD \log \frac{\gw^n}{\to^n} + g^{k \bl} g^{m \bn}
g^{p \bq} \left( g_{m \bq, k} - g_{k \bq, m} \right) \left(
g_{\bn p, \bl} - g_{\bl p, \bn} \right)\\
&\ + \tr_{\gw} \del \delb \log \det \til{g} - \tr_{\to} \psi.
\end{align*}
Lastly we compute
\begin{align*}
\dt \tr_{\gw} \to =&\ g^{i \bq} g^{p \bj} P^{1,1}_{p \bq} \til{g}_{i \bj} +
\tr_{\gw} \psi\\
=&\ g^{i \bq} g^{p \bj} g^{r \bs} \left[ - g_{p \bq, r \bs} + g^{u \bv} g_{p
\bv, r} g_{u \bq, \bs} \right] \til{g}_{i \bj}\\
&\ - g^{i \bq} g^{p \bj} g^{r \bs} g^{u \bv} \left( g_{r \bv, p} - g_{p \bv, r}
\right) \left( \del_{\bq} g_{u \bs} - \del_{\bs} g_{u \bq} \right) \til{g}_{i
\bj} + \tr_{\gw} \psi.
\end{align*}
\end{proof}
\end{lemma}

\noindent We next record a lemma fixing certain canonical coordinates for $g$.
\begin{lemma} \label{coordinates} (\cite{GL} Lemma 2.1) Fix $(M^{2n}, \til{g},
J)$ a complex
manifold with Hermitian metric, and $g$ another Hermitian metric on $M$.  For $p
\in M$, there exist complex coordinates
near $p$ such that
\begin{align*}
\til{g}_{i \bj}(p) =&\ \gd_{ij}, \qquad \del_j \til{g}_{i \bi} = 0, \qquad
{g}_{i \bj} =
{g}_{i \bi} \gd_{ij}.
\end{align*}
\end{lemma}

\begin{prop} \label{PCFtraceev} Let $(M^{2n}, \til{\omega}(t), J)$ be a compact
complex
manifold with a one-parameter family of pluriclosed metrics and suppose
$\omega(t)$ is a solution of
(\ref{PCF}), and $\phi$ is a solution to (\ref{potentialeqn}).  Fix a constant
$A > 0$ and let
\begin{align*}
F =&\ \log \tr_{\to} \gw - A \phi.
\end{align*}
There is a constant $C = C(\til{g})$ such that
\begin{align*}
\left(\gD - \dt \right) F \geq&\ - \sum_i \frac{1}{\tr_{\to} \gw}
g^{m \bn} {g}^{p \bq} \left( {g}_{m \bq, i} - {g}_{i \bq,
m} \right) \left( {g}_{\bn p, \bi} - {g}_{\bi p, \bn} \right)\\
&\ + \left( A
- C \right) \tr_{\gw} \to + \frac{\left< \psi, \omega \right>_{\to}}{\tr_{\to}
\omega} - An.
\end{align*}
\begin{proof} Fix a point $p \in M$, choose coordinates for $\til{g}$ centered
at $p$ as in
Lemma \ref{coordinates} and
compute
\begin{align*}
\gD \tr_{\to} \gw =&\ g^{i \bj} \del_i \del_{\bj} \left( \til{g}^{k \bl}
g_{k \bl} \right) = \sum_{i,k} g^{i \bi} g_{k \bk, i \bi} - 2 \Re \left(
\sum_{i,j,k}
g^{i \bi} \til{g}_{j \bk, \bi} g_{k \bj, i} \right) + \mathcal
O\left((\tr_{\omega} \til{\omega}) (\tr_{\til{\omega}} \omega) \right).
\end{align*}
Plugging this calculation into the result of Lemma \ref{prelimcalcs}, we
conclude
\begin{align*}
\left(\gD - \dt \right) \tr_{\to} \gw =&\ - 2 \Re \left(
\sum_{i,j,k} g^{i
\bi} \til{g}_{j \bk, \bi} g_{k \bj, i} \right) + \sum_k g^{i \bi} g^{m \bm}
\left( g_{k \bm, i} g_{\bk m, \bi} \right)\\
&\ - {g}^{m \bn} {g}^{p \bq} \left( {g}_{m \bq, i} - {g}_{i \bq,
m} \right) \left( {g}_{\bn p, \bi} - {g}_{\bi p, \bn} \right) + \mathcal
O\left((\tr_{\omega} \til{\omega}) (\tr_{\til{\omega}} \omega)
\right) + \left< \psi, \omega \right>_{\to}.
\end{align*}
Using the properties of Lemma \ref{coordinates} we can estimate
\begin{align*}
\brs{ -2 \Re \sum_{i,j,k} {g}^{i \bi} \til{g}_{j \bk, \bi} {g}_{k \bj, i}}
\leq&\ \sum_{i, j \neq k} {g}^{i \bi} {g}^{j \bj} {g}_{k \bj,
i} g_{\bk j, \bi} + \mathcal O((\tr_{\omega} \til{\omega})(\tr_{\to} \omega)).
\end{align*}
Next let us apply the properties of the coordinates in Lemma \ref{coordinates}
and the Cauchy-Schwarz inequality twice to yield
\begin{align*}
\frac{\brs{\del \tr_{\to} \gw}^2_{g}}{\tr_{\to} \gw} =&\ \frac{1}{\tr_{\to} \gw}
\sum_{i,j,k} {g}^{i \bi} \del_i {g}_{j \bj} \del_{\bi}
{g}_{k \bk}\\
=&\ \frac{1}{\tr_{\to} \gw} \sum_{j,k} \sum_i \sqrt{{g}^{i \bi}} \del_i
{g}_{j \bj} \sqrt{{g}^{i \bi}} \del_{\bi} {g}_{k \bk}\\
\leq&\ \frac{1}{\tr_{\to} \gw} \sum_{j,k} \left( \sum_i {g}^{i \bi}
\brs{\del_i {g}_{j \bj}}^2 \right)^{\frac{1}{2}} \left( \sum_i {g}^{i
\bi} \brs{\del_i {g}_{k \bk}}^2 \right)^{\frac{1}{2}}\\
=&\ \frac{1}{\tr_{\to} \gw} \left( \sum_j \left( \sum_i {g}^{i \bi}
\brs{\del_i {g}_{j \bj}}^2 \right)^{\frac{1}{2}} \right)^2\\
=&\ \frac{1}{\tr_{\to} \gw} \left( \sum_j \sqrt{{g}_{j\bj}} \left( \sum_i
{g}^{i \bi} {g}^{j \bj} \brs{\del_i {g}_{j \bj}}^2
\right)^{\frac{1}{2}} \right)^2\\
\leq&\ \sum_{i,j} {g}^{i \bi} {g}^{j \bj} \del_i {g}_{j \bj}
\del_{\bi} {g}_{j \bj}.
\end{align*}
Finally we can now conclude
\begin{align*}
\sum_k {g}^{i \bi} {g}^{m \bm} & \left( {g}_{k \bm, i} {g}_{\bk m, \bi}
\right) - \frac{\brs{\del \tr_{\to} {g}}_{{g}}^2}{\tr_{\to} {\gw}} - 2 \Re
\left( \sum_{i,j,k} {g}^{i
\bi} \til{g}_{j \bk, \bi} {g}_{k \bj, i} \right) \geq - \mathcal O((\tr_{\omega}
\to) (\tr_{\to} \omega)).
\end{align*}
Combining the above calculations yields the result.
\end{proof}
\end{prop}

\begin{prop} \label{tracevolumeev} Let $(M^{2n}, \til{\gw}(t), J)$ be a compact
complex
manifold with a one-parameter family of pluriclosed metrics and suppose $\gw(t)$
is a solution of
(\ref{PCF}), and $\phi$ is a solution to (\ref{potentialeqn}).  Fix a constant
$A > 0$ and let
\begin{align*}
F =&\ \log \tr_{\to} \gw - \log \frac{\gw^n}{\to^n} - A \phi.
\end{align*}
There is a constant $C = C(\til{g})$ such that
\begin{align*}
\left( \gD - \dt \right) F \geq \left( A
- C \right) \tr_{\gw} \to - \frac{\left< \psi, \omega \right>_{\to}}{\tr_{\to}
\omega} - \tr_{\to} \psi - An.
\end{align*}
\begin{proof} Starting from the result of Proposition \ref{PCFtraceev} and using
Lemma \ref{prelimcalcs} we immediately conclude
\begin{align*}
\left({\gD} - \dt \right) F \geq&\ {g}^{k \bl} {g}^{m \bn}
{g}^{p \bq} \left( {g}_{m \bq, k} - {g}_{k \bq, m} \right) \left(
{g}_{\bn p, \bl} - {g}_{\bl p, \bn} \right) \\
&\ - \frac{1}{\tr_{\to} \gw} {g}^{m \bn} {g}^{p \bq} \left(
{g}_{m \bq, i} - {g}_{i \bq,
m} \right) \left( {g}_{\bn p, \bi} - {g}_{\bi p, \bn} \right)\\
&\ + \left( A
- C \right) \tr_{\gw} \to + \frac{\left< \psi, \omega \right>_{\to}}{\tr_{\to}
\omega} - \tr_{\to} \psi - An.
\end{align*}
Let ${T}_{k m \bq} = {g}_{m \bq, k} - {g}_{k \bq, m}$.  The first
term above is then $\brs{T}^2$.  Likewise, if we choose
coordinates at a point such that $\til{g}_{i \bj} = \gd_{ij}$ and $g$ is
diagonalized, specifically $g_{i \bi} = \gl_i$ then
\begin{align*}
\frac{1}{\tr_{\to} \gw} \leq&\ \frac{1}{\gl_i}
\end{align*}
for any $i$.  Thus
\begin{align*}
- \frac{1}{\tr_{\to} \gw} {g}^{m \bn} {g}^{p \bq} \left( {g}_{m
\bq, i} - {g}_{i \bq,
m} \right) \left( {g}_{\bn p, \bi} - {g}_{\bi p, \bn} \right) =&\ -
\frac{1}{\tr_{\to} \gw} \sum_i {g}^{m \bn} {g}^{p \bq} {T}_{i m
\bq} {T}_{\bi \bn p}\\
\geq&\ - \sum_i {g}^{i \bi} {g}^{m \bn} {g}^{p \bq} {T}_{i m
\bq} {T}_{\bi \bn p}\\
=&\ - \brs{{T}}^2.
\end{align*}
The proposition follows.
\end{proof}
\end{prop}

\begin{prop} \label{inversemetricev} Let $(M^{2n}, \til{\gw}(t), J)$ be a
compact
complex
manifold with a one-parameter family of pluriclosed metrics and suppose $\gw(t)$
is a solution of
(\ref{PCF}).  Then in local complex coordinates
\begin{align*}
\left(\gD - \dt \right) \tr_{\gw} \til{\gw} =&\ g^{r \bs} g^{i \bv} g^{u \bq}
g^{p \bj} g_{u \bv, r} g_{p \bq, \bs} \til{g}_{i \bj} - 2 \Re \left(g^{r \bs}
g^{i \bq} g^{p \bj} g_{p \bq, \bs} \til{g}_{i \bj, r} \right)+ g^{r \bs} g^{i
\bj} \til{g}_{i \bj, r \bs}\\
&\ + g^{i \bq} g^{p \bj} g^{r \bs} g^{u \bv} \left( g_{r \bv, p} - g_{p \bv, r}
\right) \left( g_{u \bs,\bq} - g_{u \bq,\bs} \right) \til{g}_{i
\bj} - \tr_{\gw} \psi.
\end{align*}
\begin{proof} We directly compute
\begin{align*}
\gD \tr_{\gw} \to =&\ g^{r \bs} \del_r \del_{\bs} \left[ g^{i \bj} \til{g}_{i
\bj} \right]\\
=&\ g^{r \bs} \del_r \left[ - g^{i \bq} g^{p \bj} g_{p \bq, \bs} \til{g}_{i \bj}
+ g^{i \bj} \til{g}_{i \bj, \bs} \right]\\
=&\ g^{r \bs} \left[ g^{i \bv} g^{u \bq} g^{p \bj} g_{u \bv, r} g_{p \bq, \bs}
\til{g}_{i \bj} + g^{i \bq} g^{p \bv} g^{u \bj} g_{u \bv, r} g_{p \bq, \bs}
\til{g}_{i \bj} - g^{i \bq} g^{p \bj} g_{p \bq, r \bs} \til{g}_{i \bj} \right.\\
&\ \qquad \left. - g^{i \bq} g^{p \bj} g_{p \bq, \bs} \til{g}_{i \bj, r} - g^{i
\bq} g^{p \bj} g_{p \bq, r} \til{g}_{i \bj, \bs} + g^{i \bj} \til{g}_{i \bj, r
\bs} \right].
\end{align*}
Combining this with the result of Lemma \ref{prelimcalcs} yields the result.
\end{proof}
\end{prop}

\begin{prop} \label{Linf} Let $(M^{2n}, \til{\gw}(t), J)$ be a compact complex
manifold with a one-parameter family of 
pluriclosed
metrics and suppose $\gw(t)$ is a solution of (\ref{PCF}) on $[0, \tau)$,
and $\phi$ is a solution to (\ref{potentialeqn}).  There is a constant $C =
C(\til{g}, \psi)> 0$
such that
\begin{align*}
\tr_{\to} \gw \leq&\ C e^{ C \left(t + \log \frac{\gw^n}{\to^n} + \phi \right)}.
\end{align*}
\begin{proof} Let
\begin{align*}
 F = \log \tr_{\til{\gw}} \gw - \log \frac{\gw^n}{\til{\gw}^n} - A \phi - B t,
\end{align*}
where the constants $A$ and $B$ are to be determined.  Using Proposition
\ref{tracevolumeev} we obtain
\begin{align*}
 \left( \gD - \dt \right) F \geq \left( A
- C \right) \tr_{\gw} \to - \frac{\left< \psi, \omega \right>_{\to}}{\tr_{\to}
\omega} - \tr_{\to} \psi - An + B.
\end{align*}
Now we note that
\begin{align*}
\brs{ \left< \psi, \gw \right>_{\til{\gw}}} \leq \brs{\psi}_{\til{\gw}}
\brs{\gw}_{\til{\gw}} \leq C \brs{\psi}_{\til{\gw}} \tr_{\til{\gw}} \gw.
\end{align*}
Thus we have
\begin{align*}
 \left( \gD - \dt \right) F \geq \left(A - C - C \brs{\psi}_{\til{\gw}} \right)
- An + B.
\end{align*}
Choosing $A$ sufficiently large with respect to the constants and
$\brs{\psi}_{\til{\gw}}$ and then choosing $B = n A$ yields
\begin{align*}
 \left( \gD - \dt \right) F \geq 0.
\end{align*}
By the maximum principle, for any $t < \tau$ we conclude
\begin{align*}
\sup_M F(t) \leq \sup_M F(0).
\end{align*}
The result follows.
\end{proof}
\end{prop}

\section{\texorpdfstring{$C^1$}{C1} metric estimate}
In this section we derive certain differential inequalities for the Chern
connection along a solution to (\ref{PCF}).  These estimates are inspired by
Calabi's third-order estimate of the potential function for the complex Monge
Ampere equation, and similar estimates for the K\"ahler-Ricci flow were
considered in \cite{PSS}.  Fix $(M^{2n}, \til{g}, J)$ a Hermitian manifold and
let $g$ denote another Hermitian metric on $M$.  Let $h$ denote the endomorphism
of the
tangent bundle
\begin{align*}
h^i_j = \til{g}^{i \bk} g_{\bk j}.
\end{align*}
Let
\begin{align*}
\U = \N h h^{-1},
\end{align*}
and let
\begin{align*}
W = \brs{\U}^2 =&\ g^{i \bj} g^{k \bl} g_{m \bn}  \U_{i k}^m \bar{ \U _{j l}^n}
\end{align*}
where the first lowered index on the tensor $\U = \N h h^{-1}$ is that arising
from
the derivative.  The tensor $\U$ is the difference of the Chern connections
induced by $\til{g}$ and $g$.  In particular one observes that
\begin{align*}
\bar{\N} \U =&\ \til{\Omega} - \Omega.
\end{align*}

\begin{lemma} \label{Ulaplace} Let $(M^{2n}, \til{g}, J)$ be a Hermitian
manifold, let $g$ denote
another Hermitian metric on $M$, and let $W$ be defined as above.  Then
\begin{align*}
\gD W =&\ \brs{ \bar{\N} \U}^2 + \brs{\N \U}^2\\
&\ + g^{i \bj} g^{k \bl} g_{m \bn} \bar{\left( - g^{\bq p} T_{j p}^r \Omega_{\bq
r
l}^n - \N_j S_l^n + g^{p \bq} \N_p \til{\Omega}_{\bq j l}^n \right)} \U_{i
k}^{m} +
\mbox{ conjugate }\\
&\ + \U_{i k}^m \left( S^{i \br} g^{k \bl} g_{m \bn}
\bar{ \U_{r l}^n} + g^{i \bj} S^{k \br} g_{m \bn} \bar{
\U_{j r}^n} - g^{i \bj} g^{k \bl} S_{r \bm}
\bar{\U_{j l}^m} \right).
\end{align*}
\begin{proof}
First we compute
\begin{align*}
\gD W =&\ g^{i \bj} g^{k \bl} g_{m \bn} \left[ \gD \U_{i
k}^m \bar{ \U_{j l}^n} + \U_{i k}^m 
\bar{ \bar{\gD} \U_{j l}^n} \right] + \brs{ \bar{\N} \U}^2 + \brs{\N \U}^2.
\end{align*}
Next we commute derivatives to see
\begin{align*}
\bar{\gD} \U_{j l}^n =&\ g^{p \bq} \N_{\bq} \N_{p}
\U_{j l}^n\\
=&\ g^{p \bq} \left[ \N_{p} \N_{\bq} \U_{j l}^n -
\Omega_{\bq p j}^r \U_{r l}^n - \Omega_{\bq p l}^r
\U_{j r}^n + \Omega_{\bq p r}^n \U_{j l}^r \right]\\
=&\ \gD \U_{j l}^n + S_j^r \U_{r
l}^n + S_{l}^r \U_{j r}^n - S_{r}^n \U_{j l}^r
\end{align*}
Finally we observe using a general formula for curvatures of Hermitian metrics
that
\begin{align*}
\gD \U_{j l}^n =&\ g^{\bq p} \N_p \N_{\bq} \U_{j l}^n\\
=&\ g^{\bq p} \N_p \left( \til{\Omega}_{\bq j l}^n - {\Omega}_{\bq j l}^n
\right)\\
=&\ g^{\bq p} \left( - \N_j \Omega_{\bq p l}^n - T_{j p}^r \Omega_{\bq r l}^n +
\N_{p} \til{\Omega}_{\bq j l}^n \right)\\
=&\ - \N_j S_l^n - g^{\bq p} T_{j p}^r \Omega_{\bq r l}^n + g^{\bq p} \N_p 
\til{\Omega}_{\bq j l}^n.
\end{align*}
Combining these calculations yields the result.
\end{proof}
\end{lemma}

\begin{prop} \label{Wevolution} Let $(M^{2n}, \til{g}, J)$ be a Hermitian
manifold, let $g$ denote
a solution to (\ref{PCF}), and let $W$ be defined as above.  Then
\begin{align*}
\left(\gD - \dt \right) W =&\ \brs{ \bar{\N} \U}^2 +
\brs{\N \U }^2\\
&\ + g^{i \bj} g^{k \bl} g_{m \bn} \bar{\left( - g^{\bq p} T_{j p}^r \Omega_{\bq
r
l}^n - \N_j Q_l^n + g^{p \bq} \N_p \til{\Omega}_{\bq j l}^n \right)} \U_{i
k}^{m}\\
&\ + \mbox{ \emph{conjugate} }\\
&\ + \U_{i k}^m \left( (Q^1)^{i \br} g^{k \bl} g_{m \bn}
\bar{ \U_{r l}^n} + g^{i \bj} (Q^1)^{k \br} g_{m \bn} \bar{
\U_{j r}^n} - g^{i \bj} g^{k \bl} Q^1_{r \bm} \bar{\U_{j l}^m} \right)
\end{align*}
\begin{proof}
First observe the variational equation
\begin{align*}
\dt \left( \N h h^{-1} \right) =&\ \N \left( h^{-1} \left( \dt h \right)
\right).
\end{align*}
It follows from Lemma \ref{Ulaplace} that for a general variation one has
\begin{align*}
\left(\gD - \dt \right) W =&\ \brs{ \bar{\N} \U}^2 +
\brs{\N \U}^2\\
&\ + g^{i \bj} g^{k \bl} g_{m \bn} \bar{\left( - g^{\bq p} T_{j p}^r \Omega_{\bq
r
l}^n - \N_j S_l^n - \N_j \left( h^{-1} \dot{h} \right)_l^n + g^{p \bq} \N_p
\til{\Omega}_{\bq j l}^n \right)} \U_{i k}^{m}\\
&\ + \mbox{ conjugate }\\
&\ + \U_{i k}^m \left( \left( h^{-1} \dot{h}^{i \br} +
S^{i \br} \right) g^{k \bl} g_{m \bn} \bar{ \U_{r l}^n}
\right.\\
&\ \qquad \left.  + g^{i \bj} \left( \left( h^{-1} \dot{h} \right)^{k \br} +
S^{k \br} \right) g_{m \bn} \bar{ \U_{j r}^n} - g^{i \bj} g^{k \bl} \left(
h^{-1} \dot{h}_{r \bm} + S_{r \bm}
\right) \bar{\U_{j l}^m} \right)
\end{align*}
Plugging in $\dot{g} = -P^{1,1} = - S + Q^1$ yields the result.
\end{proof}
\end{prop}

\begin{prop} \label{Westimate} Let $(M^{2n}, \til{g}, J)$ be a Hermitian
manifold, let $g$ denote
a solution to (\ref{PCF}) on $[0, \tau)$, and let $W$ be defined as above. 
Suppose there exists a constant $K$ such that
\begin{align*}
\frac{1}{K} \til{g} \leq g(t) \leq K \til{g}
\end{align*}
for all $ t \in [0, \tau)$.  Then there is a constant $C(K, \til{g})$ such that
\begin{align*}
\left(\gD - \dt \right) W \geq&\ - C \left(1 + W  + \brs{T}^2 W \right).
\end{align*}

\begin{proof} We start by noting that
\begin{align*}
\bar{\N} \U =&\ \Omega - \til{\Omega}, 
\end{align*}
thus
\begin{align*}
\brs{\bar{\N} \U}^2 \geq&\ \frac{1}{2} \brs{\Omega}^2 - C \brs{\til{\Omega}}^2.
\end{align*}
Also, by orthogonally projecting $\U$ onto its skew symmetric part, one observes
that
\begin{align*}
\brs{\bar{\N} \U}^2 + \brs{\N \U}^2 \geq&\ \frac{1}{2} \brs{\bar{\N} T}^2 +
\frac{1}{2} \brs{\N T}^2 - C \brs{\N \til{T}}^2\\
\geq&\ \frac{1}{2} \brs{\bar{\N} T}^2 +
\frac{1}{2} \brs{\N T}^2 - C W.
\end{align*}
Thus, starting from the result of Proposition \ref{Wevolution} we first estimate
using the Cauchy-Schwarz inequality
\begin{align*}
g^{i \bj} g^{k \bl} g_{m \bn} g^{p \bq} \N_p \til{\Omega}_{\bq j l}^n \U_{i k}^m
=&\ g * g^{-1} \left( \til{\N} \til{\Omega} + \U * \til{\Omega} \right) * \U\\
\leq&\ C(K) \left[ \brs{\til{\N} \til{\Omega}}^2_{g} + W \right]\\
\leq&\ C(K, \til{g}) \left(1 + W \right).
\end{align*}
Next we estimate
\begin{align*}
\brs{g^{i \bj} g^{k \bl} g_{m \bn} g^{p \bq} T_{j p}^r \Omega_{\bq r l}^n \U_{i
k}^m} \leq&\ C \brs{T} \brs{\Omega} \brs{\U}\\
\leq&\ \theta \brs{\Omega}^2 +\frac{C}{\theta} \brs{T}^2 \brs{\U}^2\\
\leq&\ \theta \brs{\Omega}^2 + \frac{C}{\theta} \brs{T}^2 W.
\end{align*}
Similarly we have
\begin{align*}
\brs{ g^{i \bj} g^{k \bl} g_{m \bn} \N_j Q_l^n \U_{i k}^n} \leq&\ C
\left[\brs{\N T} + \brs{\bar{\N} T} \right] \brs{T} \brs{\U}\\
\leq&\ \theta \left[ \brs{\bar{\N} T}^2 + \brs{\N T}^2 \right] +
\frac{C}{\theta} \brs{T}^2 W.
\end{align*}
Finally, one clearly has
\begin{align*}
\U_{i k}^m Q^{i \br} g^{k \bl} g_{m \bn} \bar{\U_{r l}^n} \leq&\ C\brs{T}^2W,
\end{align*}
and likewise for the rest of the terms.  Choosing $\theta$ sufficiently small
and combining these estimates yields the result.
\end{proof}
\end{prop}

\section{Regularity Theorems}

In this section we give the proofs of the regularity theorems stated in the
introduction.  We start by giving the proof of Theorem \ref{ricciregularity}.

\begin{proof}
\begin{lemma} \label{PCFregularitylemma1} Let $(M^{2n}, \omega(t), J)$ be a
solution to (\ref{PCF}) on a finite time interval $[0, \tau)$
and suppose 
\begin{align*}
\int_0^\tau \brs{\frac{\del g}{\del t}} < \infty.
\end{align*}
Then there is a $C^0$ metric $\omega(\tau)$ such that
\begin{align*}
\lim_{t \rightarrow \tau} \omega(t) =&\ \omega(\tau)\\
\end{align*}
in $C^0$.
\begin{proof} This is Lemma 14.2 of \cite{Ham3}.
\end{proof}
\end{lemma}

\begin{prop} \label{PCFregularityprop1} Let $(M^{2n}, \omega(t), J)$ be a
solution to (\ref{PCF}) on $[0,
\tau)$.  Suppose $\omega(t) \rightarrow \omega(\tau)$ in $C^0$, i.e. the limit
at time $\tau$ exists as a $C^0$ metric.  Then in fact $\omega(t)$ is bounded in
$C^1$ and for all $p < \infty$ there are constants $C_p$ such that
\begin{align*}
\int_M \brs{\Omega(\omega_\tau)}^p \leq&\ C_p
\end{align*}
\begin{proof}  For the $C^1$ norm on metrics we choose finitely many local
coordinate patches to cover $M$ and take the supremum over the coordinate
derivatives in all these charts.  It is equivalent to choose a fixed connection
$\N_0$ and use the covariant derivative with respect to that fixed connection.

Suppose that $\{\omega(t)\}$ is unbounded in $C^1$.  Let $\phi(x,t) =
\brs{\N_0 \omega}$.  Then for some sequence
$(x_i, t_i)$, where $t_i \rightarrow \tau$, we have $\sup_{M} \phi(t_i)$ is
achieved at $x_i$, and moreover goes to infinity.  By choosing
a subsequence we may assume that $x_i$ converges to some point $x \in M$.  By
choosing a coordinate chart around $x$, and translating coordinates for $i$
sufficiently large we may assume that $sup_{z} \phi = \ga_i$ is
attained at $z = 0$.  Now choose new coordinates $w = \ga_i z$.  Since $\omega$
is converging in $C^0$, it follows that
\begin{align} \label{C0proofloc10}
\lim_{i \rightarrow \infty} \int_{\{\brs{w} < 1\}} \phi = 0.
\end{align}
One may express the system (\ref{PCF}) in local coordinates as
\begin{gather} \label{PCFcoordinates}
\begin{split}
\dt \omega_{ij} - g^{k l} \del_k \del_l \omega_{ij}
=&\ \omega^{-1} * \omega^{-1} \left( \del \omega^{*2} \right).
\end{split}
\end{gather}

The right hand side of the equation is uniformly bounded in $C^0$ on
$[0, \tau)$, hence each coordinate
function is the solution to a uniformly parabolic equation with continuous
coefficients and bounded right hand side.  It follows by \cite{Lieberman}
Theorem 7.13 that on $\{|w| < 1 - \ge \}$ we have $\brs{\omega}_{H_2^p} <
\infty$ for all $p$.

Choosing $p > 2n$, and applying the Sobolev inequality, we
attain a uniform $C^1$ bound for $\omega(t_i)$, and moreover a
convergent
subsequence at time $\tau$.  But then, for this subsequence,
(\ref{C0proofloc10}) implies that $\lim_{i \rightarrow \infty} \phi(t_i)_{w = 0}
= 0$, a contradiction.  It follows that $\omega(t)$ is
bounded in $C^1$.  Now applying the above regularity argument to $\omega$ in
local coordinates around any point yields again the $H_2^p$ bound on $\omega$
for all $p$, and hence the curvature $\Omega$ is bounded in $L^p$ for all $p$.
\end{proof}
\end{prop}

We now proceed with the main argument.  Suppose that the statement of the
theorem were false.  Then let $(M^{2n},
\omega(t), J)$ be a solution to (\ref{PCF}) on $[0, \tau)$
satisfying
\begin{align*}
\int_0^{\tau} \brs{\frac{\del g}{\del t}} < \infty.
\end{align*}
By Lemma \ref{PCFregularitylemma1} we conclude the existence of a $C^0$ limit
metric $\omega(\tau)$.  By Proposition \ref{PCFregularityprop1} we conclude that
in fact $\omega(t)$ is a $C^1$ metric and furthermore one has uniform $L^p$
bounds on the curvature as $t \rightarrow \tau$.  By (\cite{ST1} Theorem 1.1),
if
we can show a uniform bound on the $C^0$ norm of curvature and torsion, we can
conclude smooth convergence of the metrics as $t \rightarrow \tau$.  By the
general short-time existence result for these equations, we conclude
that $\tau$ is not the maximal existence time, providing the contradiction.

So, we can differentiate (\ref{PCFcoordinates}) to yield
\begin{align*}
\dt \left( \N_0 \omega \right)_{ijk} - g^{pq} \del_p \del_q (\N_0 \omega)_{ijk}
=&\ \omega^{-1} * \omega^{-1} * \left( \del \omega * \del^2 \omega \right)
\end{align*}
By the discussion above, the right hand side is uniformly bounded in $L^p$, so
again we conclude that $\omega$ is uniformly bounded in $H_3^p$.  Choosing $p$
sufficiently large and applying the Sobolev inequality, we conclude a uniform
$C^0$ bound on $\Omega$, $T$ and $\N T$, and the theorem follows.
\end{proof}

\noindent Now we give the proof of Theorem \ref{regularity2}.

\begin{proof} Consider the differential inequality of Proposition
\ref{PCFtraceev}.  First we note that
\begin{align*}
\frac{1}{\tr_{\to} \omega} g^{m \bn} g^{p \bq} \left( g_{m \bq, i} - g_{i \bq,
m} \right) \left( g_{\bn p, \bi} - g_{\bi p, \bn} \right) \leq&\ \brs{T}^2.
\end{align*}
Applying the maximum principle, since $\phi$ and $\brs{T}^2$ are bounded we
conclude a uniform upper bound on $\tr_{\to} \gw$.  In particular, there is a
constant $K$ such that
\begin{align*}
g(t) \leq K \til{g}.
\end{align*}
Also, by using Lemma \ref{prelimcalcs}, if $F = \log \frac{\det g}{\det \til{g}}
+ A \phi$, we compute that
\begin{align*}
 \dt F \geq&\ \gD F + \brs{T}^2 + \left(A - C\right) \tr_{\omega} \til{\omega} -
n.
\end{align*}
In particular, for $A$ chosen sufficiently large we conclude by the maximum
principle,
\begin{align*}
\inf_{M \times [0, \tau)} \log \frac{\det g}{\det \til{g}} \geq - C(\brs{\phi}).
\end{align*}
We thus conclude a uniform lower bound for $g(t)$ on $[0, \tau)$ by the
arithmetic-geometric mean.  Again using that the torsion is bounded, we may
apply
the maximum principle to the differential inequality in Proposition
\ref{Westimate} to conclude a uniform $C^1$ bound on $g(t)$ on $[0, \tau)$. 
Equation (\ref{PCF}) is strictly parabolic in local complex coordinates, with
bounded $C^1$ norm, so uniform $C^k$ estimates follow from the Schauder theory
for all $k$, and the theorem follows.
\end{proof}

Let us finish this section with a few remarks on the nature of the potential
function $\phi$.  First, Theorem \ref{regularity2} even provides a slightly
different perspective on the regularity of K\"ahler Ricci flow.  In this case
the torsion vanishes for all time, so one only has to check that the potential
function is bounded.  It is clear by applying the maximum principle to
(\ref{potentialeqn}) that if
\begin{align*}
\int_0^{\tau} \sup_M \tr_{\gw} \til{\gw} < \infty,
\end{align*}
then the flow will extend past time $\tau$.  This condition can be checked in
certain settings.  Also, in the general non-K\"ahler setting the function $\phi$
is automatically bounded on certain background manifolds.  Consider the
following lemma.

\begin{lemma} Let $(M^{2n}, \til{\omega}, J)$ be a complex manifold and suppose
$\til{\omega}$ is K\"ahler and moreover $\til{\Omega} \leq 0$, in the sense of
sections of $\End (\Lambda^{1,1})$.  Let $\omega(t)$ denote a solution to
(\ref{PCF}) on $[0, \tau)$.  Then there is a constant $C > 0$ so that
\begin{align*}
\tr_{\gw} \to \leq&\ C\\
\brs{\phi} \leq C\tau.
\end{align*}
\begin{proof} From Proposition \ref{inversemetricev}, if $\til{\omega}$ is
K\"ahler we can choose coordinates where $\del_i\til{g}_{j \bk} = 0$ and
simplify, since $Q^1 \geq 0$,
\begin{align*}
\left(\gD - \dt \right) \tr_{\gw} \til{\gw} \geq&\ - g^{r \bs} g^{i \bj}
\til{\Omega}_{i \bj r \bs}\\
\geq&\ 0.
\end{align*}
Applying the maximum principle proves the uniform upper bound for $\tr_{\gw}
\to$, and then the bound for $\phi$ follows immediately applying the maximum
principle to (\ref{potentialeqn}).
\end{proof}
\end{lemma}

Since $\til{\omega}$ is K\"ahler, $P$ is just the Ricci form, hence the
hypotheses are satisfied on complex tori or K\"ahler manifolds with negative
curvature operator.  Note of course that we are not assuming $\omega$ is
K\"ahler.  This bound suggests that the function $\phi$ is quite natural to
introduce, and furthermore suggests that its possible blowup is not related to
any ``local'' singularity model since it is bounded on these natural background
manifolds.

\section{Conjectural Picture of Singularity Formation} \label{conesection}

The notion of the K\"ahler cone in $H^2(M, \mathbb R) \cap H^{1,1}(M, \mathbb
C)$ is crucial to understanding the structure of solutions of K\"ahler Ricci
flow.  Recall from the introduction that along a solution to K\"ahler-Ricci flow
the K\"ahler class
satisfies an ODE, depending on the normalization.  Clearly a necessary
condition for existence of the flow is that this ODE stay in the K\"ahler cone. 
As mentioned in the introduction, this condition is in fact sufficient
(\cite{TZ} Proposition 1.1).  It is natural
to conjecture that a similar phenomenon is at play guiding the singular behavior
of solutions of (\ref{PCF}).

First of all recall from the introduction that
\begin{align*}
\mathcal H^{1,1}_{\del + \delb} = \frac{ \{ \Ker \del \delb : \Lambda_{\mathbb
R}^{1,1} \rightarrow \Lambda_{\mathbb R}^{2,2} \} }{ \{ \del \ga + \delb
\bar{\ga} | \ga \in \Lambda^{0,1} \} }.
\end{align*}
This is known as the $(1,1)$ Aeppli cohomology group and one basic fact is that
this space is finite dimensional, as can be seen by constructing the necessary
short exact sequence of coherent sheaves.  Let the \emph{positive cone} inside
$\mathcal H^{1,1}_{\del +
\delb}$ be
\begin{align*}
\mathcal P_{\del + \delb} = \{ [\phi] \in \mathcal H^{1,1}_{\del + \delb} |
\exists \psi \in [\phi], \psi > 0  \}.
\end{align*}
It is clear that a necessary condition for a solution to (\ref{PCF}) to exist is
that the class $[\omega_t] = [\omega_0 - t c_1] \in \mathcal P_{\del + \delb}$. 
We state this for emphasis.

\begin{prop} Let $(M^{2n}, g_0, J)$ be a compact complex manifold with
pluriclosed metric.  Let
\begin{align*}
\tau^* := \sup_{t \geq 0} \{t | [\omega_0 - t c_1] \in \mathcal P_{\del + \delb}
\}.
\end{align*}
Let $\tau$ denote the maximal existence time of the solution of (\ref{PCF}) with
initial condition $g_0$.  Then
\begin{align*}
\tau \leq \tau^*.
\end{align*}
\end{prop}

\noindent Furthermore, in analogy with K\"ahler-Ricci flow, it is natural to
conjecture
that membership in this cone suffices for existence.

\begin{conj} \label{existenceconj} \textbf{\emph{Weak existence conjecture:}}
Let $(M^{2n}, g_0, J)$ be a compact complex
manifold with
pluriclosed metric.  Let
\begin{align*}
\tau^* := \sup_{t \geq 0} \{ t | [\omega_0 - t c_1] \in \mathcal P_{\del +
\delb} 
\}.
\end{align*}
Then the solution to (\ref{PCF}) with initial condition $g_0$ exists on
$[0, \tau^*)$, and $\tau^*$ is the maximal time of existence.
\end{conj}

Let us note here that this we are implicitly making this conjecture, and the two
related conjectures below, for any normalization of (\ref{PCF}), i.e. the volume
normalized version of (\ref{PCF}) or other possible normalizations.  A stronger
version of this conjecture would be that there are uniform
$C^{\infty}$ estimates on $\omega(t)$ depending on $d([\omega(t)], \del \mathcal
P_{\del + \delb})$, where this means distance with respect to some metric
defined on $\mathcal H^{1,1}_{\del + \delb}$.  Let us also state this for
emphasis.

\begin{conj} \label{strongexistenceconj} \textbf{\emph{Strong existence
conjecture:}} Let $(M^{2n}, g_0, J)$ be a compact complex
manifold with
pluriclosed metric.  Let $\omega(t)$ be the solution to (\ref{PCF}) with initial
condition $\omega_0$.  Let $t$ be such that $[\omega_0 - t c_1] \in \mathcal
P_{\del + \delb}$.  Then there exist uniform $C^{\infty}$ estimates on
$\omega(t)$ depending on $\omega_0, t$ and $d([\omega(t)], \del \mathcal P_{\del
+
\delb})$.  Moreover, there exist uniform bounds on the curvature and diameter of
$\gw(t)$.
\end{conj}

It is possible to characterize $\mathcal P_{\del + \delb}$ using more calculable
cohomological criteria in the case of non-K\"ahler complex surfaces, which will
allow us to derive some consequences of conjecture \ref{strongexistenceconj}. 
Our Theorem
\ref{positivityintro} follows by combining the positivity
result
of Buchdahl \cite{Buchdahl} on non-K\"ahler surfaces (see also Lamari
\cite{Lamari}), and further related results of Teleman \cite{Telemancone}.  Let
us start by stating the main theorem of
\cite{Buchdahl2}, which represents the main technical difficulty of Theorem
\ref{positivityintro}.

\begin{thm} \label{Buchdahl} (\cite{Buchdahl2} Main Theorem) Let $(M^{4},
\omega, J)$ be a complex
surface with pluriclosed metric $\omega$.  Suppose $\phi \in \Lambda^{1,1}$ is
pluriclosed and satisfies
\begin{itemize}
 \item {$\int_M \phi \wedge \phi > 0$}
\item{ $\int_M \phi \wedge \gw > 0$}
\item{ $\int_D \phi > 0$ for every irreducible effective divisor with $D \cdot D
< 0$.}
\end{itemize}
Then there exists $f \in C^{\infty}(M)$ such that $\phi + \sqrt{-1} \del \delb f
> 0$.
\end{thm}

\noindent For the statement of Theorem \ref{positivityintro} we need some
further
background. 
Recall the Bott-Chern cohomology group
\begin{align*}
H^{1,1}_{\text{BC}} = \frac{ \{\Ker d : \Lambda^{1,1}_{\mathbb R} \rightarrow
\Lambda^3_{\mathbb R} \}}{i \del \delb \Lambda^0_{\mathbb R}}.
\end{align*}
Also, define the groups
\begin{align*}
B^{1,1}_{\mathbb R} =&\ d \{ \Lambda^1_{\mathbb R} \} \cap
\Lambda^{1,1}_{\mathbb R},\\
H^{1,1}_{\mathbb R} =&\ \frac{ \{\Ker d : \Lambda^{1,1}_{\mathbb R} \rightarrow
\Lambda^3_{\mathbb R} \}}{B^{1,1}_{\mathbb R}}.
\end{align*}

\begin{lemma} Let $(M^4, \omega, J)$ be a complex
surface with pluriclosed metric.  Then there are exact sequences
\begin{align*}
0 &\rightarrow \frac{B^{1,1}_{\mathbb R}}{i \del \delb \Lambda^0_{\mathbb R}}
\rightarrow H^{1,1}_{\text{BC}} \rightarrow H^{1,1}_{\mathbb R} \rightarrow 0\\
0 &\rightarrow i \del \delb \Lambda^0_{\mathbb R} \rightarrow B^{1,1}_{\mathbb
R} \rightarrow \mathbb R,
\end{align*}
where the final map above is given by the $L^2$ inner product with $\omega$.
\begin{proof} We include the elementary proof for convenience.  The first exact
sequence of $(1)$ is tautological.  For the second sequence, fix $\mu \in
B^{1,1}_{\mathbb R}$ satisfying
\begin{align*}
\int_M \mu \wedge \omega = 0.
\end{align*}
It follows from the maximum principle that the adjoint of $\tr_{\omega} \del
\delb$ has kernel only constant functions since the adjoint operator annihilates
constants.  Thus by standard theory we can now solve
\begin{align*}
\gD u = \tr_{\omega} \mu.
\end{align*}
Thus $i \del \delb u - \mu$ is exact, and also anti-self-dual since its inner
product with $\omega$ vanishes.  Thus it vanishes, and the lemma follows.
\end{proof}
\end{lemma}

Furthermore, (see \cite{Telemancone} Lemma 2.3), if $b_1(M)$ is odd, the space
\begin{align*}
\Gamma = \frac{B^{1,1}_{\mathbb R}}{i \del \delb \Lambda^0_{\mathbb R}}
\end{align*}
is identified with $\mathbb R$ via the $L^2$ inner product with $\omega$.  Let
$\gamma_0$ denote a positive generator of $\Gamma$.  Since the space of
pluriclosed metrics on $M$ is connected, this orientation of $\Gamma$ is
well-defined.

\begin{thm} \label{positivityintro} Let $(M^4, J)$ be a
complex non-K\"ahler surface.
 Suppose $\phi \in
\Lambda^{1,1}$ is pluriclosed.  Then $[\phi] \in \mathcal P_{\del + \delb}$ if
and
only if
\begin{itemize}
\item{$\int_M \phi \wedge \gamma_0 > 0$}\\
\item{$\int_D \phi > 0 \mbox{ for every effective divisor with negative self
intersection}$.}
\end{itemize}
\begin{proof}  Suppose $\phi \notin \mathcal P_{\del + \delb}$.  Since the image
of $\del + \delb$ is closed in $\Lambda^{1,1}_{\mathbb R} \otimes L^2(M)$
(\cite{Buchdahl} Lemma 1), one
may apply the Hahn-Banach theorem to conclude the existence of a positive
closed current $P$ such that $P(\phi) \leq 0$.  We claim that the current $P$ is
represented by a convex linear combination of
$[\gamma_0]$ and irreducible effective divisors of negative self-intersection. 
This is \cite{Telemancone}
Corollary 3.6, and we include a sketch of the proof for convenience.  First we
note that using arguments from complex analysis one can show that the set of
irreducible effective divisors of negative self-intersection is finite
(\cite{Telemancone} Remark 3.3).  Let $\mathcal C$ denote the cone generated by
$\gamma_0$ and this finite set in $H^{1,1}_{BC}$.  If $P \notin \mathcal C$,
there exists a linear
hyperplane separating $P$ from $\mathcal C$.  Specifically we can find an
element of the dual space, represented by pairing against a pluriclosed form
$\psi$, such that
\begin{align*}
 \int_M \psi \wedge \gamma_0 >&\ 0, \qquad \int_D \psi > 0, \qquad P(\psi) =
\int_M \psi \wedge P < 0.
\end{align*}
One can show by direct inspection that $\psi + t \gamma_0$ satisfies the
criteria of Theorem \ref{Buchdahl}, hence there is $f$ such that $\psi + i \del
\delb f > 0$, and so since $P$ is positive, 
$P(\psi) \geq 0$, a contradiction.
\end{proof}
\end{thm}

\noindent The following proposition shows that Conjecture \ref{existenceconj}
implies
long-time existence of solutions to (\ref{PCF}) on minimal Class $\seven$
surfaces.
\begin{prop} Let $(M^4, \omega_0, J)$ be a minimal Class $\seven$ surface with
pluriclosed metric.  Then for all $t \geq 0$,
\begin{align*}
[\omega_0 - t c_1(\rho)] \in \mathcal P_{\del + \delb}.
\end{align*}
\begin{proof}  By the above theorem it suffices to show the integral
inequalities of Theorem \ref{positivityintro} for $\omega - t c_1(\rho)$, $t$
arbitrary.  Since $\gamma_0$ is exact the first
inequality is trivial.  Also, for any effective divisor $D$ we know that $\int_M
c_1(D) \wedge c_1 \leq 0$, hence
\begin{align*}
\int_D \omega - t c_1 \geq \int_D \omega > 0.
\end{align*}
The result follows.
\end{proof}
\end{prop}

Let us furthermore describe how we expect the presence of rational curves to
enter into the singularity formation of solutions to (\ref{PCF}) on Class
$\seven$ surfaces.  As can be seen by elementary calculations of the evolution
of the degree, one has that solutions to (\ref{PCF}) on Class $\seven$ surfaces
have volume growing at least quadratically in time.  Furthermore, the area of
divisors will grow as $K \cdot D$ where $K$ is the canonical class.  This
pairing is always nonnegative, and is zero on rational curves.   Thus, if we
renormalize to fix the volume, the boundary of the cone $\mathcal P_{\del +
\delb}$ should be reached by collapsing a curve.  This is made clearer
in section \ref{nonsingularsection} where we examine nonsingular solutions.  One
can observe at this
point though that according to our characterization of $\mathcal P_{\del +
\delb}$ in Theorem \ref{positivityintro}, the boundary may be reached, after
volume
normalizing, by having $\lim_{t \rightarrow \infty} \int_M \omega(t) \wedge
\gamma_0 =
0$.  In other words, perhaps it is this condition which fails, and not the
presence of a curve satisfying $K \cdot D = 0$.  The following proposition
effectively negates this possibility.

\begin{prop} \label{positivityprop} Let $(M^4, J)$ be a complex non-K\"ahler
surface.
 Suppose $\phi \in
\Lambda^{1,1}$ is pluriclosed, and satisfies
\begin{itemize}
\item{$\int_M \phi \wedge \phi > 0$}\\
\item{$\int_M \phi \wedge \gamma_0 \geq 0$}\\
\item{$\int_D \phi > 0 \mbox{ for every effective divisor with negative self
intersection}$.}
\end{itemize}
Then $\phi \in \mathcal P_{\del + \delb}$.
\begin{proof} Fix $\til{\omega}$ a pluriclosed metric on $M$.  Note that $\phi
\in \mathcal P_{\del + \delb}$ if and only if $\psi := \phi + a \gamma_0 \in
\mathcal P_{\del + \delb}$.  Now observe that, for $a > 0$,
\begin{align*}
\int_M \psi \wedge \psi = \int_M \phi \wedge \phi + 2 a \int_M \phi \wedge
\gamma_0 > 0.
\end{align*}
Also, since $\int_M \til{\omega} \wedge \gamma_0 > 0$, we may choose $a$ large
enough so that
\begin{align*}
\int_M \psi \wedge \til{\omega} =&\ \int_M \phi \wedge \til{\omega} + a \int_M
\gamma_0 \wedge \til{\omega} > 0.
\end{align*}
Finally, since $\gamma_0$ is $d$-exact,
\begin{align*}
\int_D \psi = \int_D \phi > 0.
\end{align*}
Therefore $\psi$ satisfies all three conditions of Buchdahl's positivity
criterion (\cite{Buchdahl2} pg. 1533), so in fact there is a function $f$ such
that
\begin{align*}
\psi + i \del \delb f > 0.
\end{align*}
The proposition follows.
\end{proof}
\end{prop}

In light of this proposition, we make a final conjecture, specializing
Conjecture \ref{strongexistenceconj} to the case $n = 2$.  Note that every
pluriclosed metric satisfies $\int_M \phi \wedge \gamma_0 > 0$, so the second
condition of Proposition \ref{positivityprop} is automatically satisfied at any
potentially singular time for a solution to (\ref{PCF}).
\begin{conj} \label{strongexistenceconjn2} \textbf{\emph{Strong existence
conjecture for surfaces:}} Let $(M^{4}, g_0, J)$ be a compact complex
surface with
pluriclosed metric.  Let $\omega(t)$ be the solution to (\ref{PCF}) with initial
condition $\omega_0$.  Suppose $\omega(t)$ exists on $[0, \tau)$ and that
\begin{itemize}
\item{$\displaystyle\lim_{t \rightarrow \tau} \int_M \gw \wedge \gw > 0$}\\
\item{There exists $A > 0$ so that $\frac{1}{A} < \displaystyle\lim_{t
\rightarrow \tau} \int_D \omega < A$
for every effective divisor with negative self
intersection.}
\end{itemize}
Then there exists a uniform bound on the curvature of $\gw(t)$
depending on $A$, and moreover the curvature remains bounded after the diameter
is rescaled to unit size.  

Furthermore, suppose $\gw(t)$ is the solution to
volume-normalized pluriclosed flow with initial condition $\gw_0$.  Suppose
$\omega(t)$ exists on $[0,\tau)$ and that there exists $A > 0$ so that
$\displaystyle \lim_{t \to \tau} \frac{1}{A} < \int_D \omega < A$
for every effective divisor with negative self
intersection.  Then there exists a uniform bound on the curvature of $\gw(t)$
depending on $A$, and moreover the curvature remains bounded after the diameter
is rescaled to unit size.  
\end{conj}

Next we want to exploit this cohomology picture to reduce solutions to
(\ref{PCF}) to an equation on a one-form.  We show that under certain
cohomological conditions related to
the Frolicher spectral sequence, solutions to equation (\ref{MAeqn})
automatically reduce to solutions to a certain equation on one-forms defined
below.  These cohomological conditions are automatically satisfied in the case
of complex surfaces.  We start with some preliminary lemmas.

\begin{lemma} \label{cohomologylemma} Let $(M^4, J)$ be a complex surface.  Then
the map
\begin{align*}
\del : H^1(\Omega^1) \ra H^1(\Omega^2)
\end{align*}
is the zero map.
\begin{proof} This argument is adapted from arguments in (\cite{BPV} IV Section
2). Let $\mathcal S$ denote the sheaf of closed holomorphic $1$-forms on $M$. 
There is an exact sequence of sheaves
\begin{align*}
0 \ra \mathbb C_M \ra \mathcal O_M \overset{d}{\ra} \mathcal S \ra 0.
\end{align*}
Since holomorphic forms on complex
surfaces are closed, we yield an exact sequence of cohomology groups
\begin{align*}
0 \ra H^0(\Omega^1) \ra H^1(M, \mathbb C) \ra H^1(\mathcal O_M)
\overset{\del}{\ra} H^1(\Omega^1).
\end{align*}
It follows from the signature theorem and the Riemann Roch formula (see
\cite{BPV} Theorem IV 2.7) that $b_1 = h^{1,0} + h^{0,1}$.  Therefore the third
map above is surjective, and hence the last map is the zero map.  Applying
Stokes theorem and Serre duality we can conclude that $\del : H^1(\Omega^1)
\ra H^1(\Omega^2)$ is also the zero map.
\end{proof}
\end{lemma}

\begin{lemma} \label{preliminarylemma} Let $(M^4, J)$ be a complex surface and
suppose $h^{0,2} = 0$.  Let $\ga \in \Lambda^{2,1}$ be a $d$-exact $(2,1)$-form.
 Then there
exists $\gg \in \Lambda^{2,0}$ such that $\ga = \delb \gg$.
\begin{proof} Since by assumption $\ga$ is exact, we conclude that
\begin{align*}
\ga =&\ d \gb = d \left( \gb^{2,0} + \gb^{1,1} + \gb^{0,2} \right)\\
=&\ \delb \gb^{2,0} + \del \gb^{1,1} + \delb \gb^{1,1} + \del \gb^{0,2}
\end{align*}
where the remaining terms vanish for dimensional reasons.  Decomposing this
equation into types yields the two equations
\begin{align*}
\ga =&\ \delb \gb^{2,0} + \del \gb^{1,1}\\
0 =&\ \delb \gb^{1,1} + \del \gb^{0,2}.
\end{align*}
Note that $\gb^{0,2}$ defines a class in $H^{0,2}(M)$.  Since $h^{0,2} = 0$,
there exists
$\mu^{0,1}$ such that $\gb^{0,2} = \delb \mu$.  We therefore conclude that
\begin{align*}
0 =&\ \delb \left( \gb^{1,1} - \del \mu^{0,1} \right).
\end{align*}
Therefore $\gb^{1,1} - \del \mu^{0,1}$ defines a class in $H^{1,1}(M)$.  Now, by
Lemma \ref{cohomologylemma} we know that $\del : H^1(\Omega^1) \ra
H^1(\Omega^2)$ is the zero map.  Therefore $\del \left(\gb^{1,1} - \del
\mu^{0,1} \right) = \del \gb^{1,1}$ represents the zero class in $H^1(\Omega^2)
\cong H^{2,1}(M)$.  Therefore there exists $\rho^{2,0}$ such that $\del
\gb^{1,1} = \delb \rho^{2,0}$.  Plugging this back into the above equation
yields
\begin{align*}
\ga =&\ \delb \gb^{2,0} + \delb \rho^{2,0}
\end{align*}
and the result follows.
\end{proof}
\end{lemma}

\begin{lemma} \label{reductionlemma} Let $(M^4,
J)$ be a complex surface.  Then
\begin{align*}
\Ker \{ \delb : \Lambda^{2,0} \rightarrow \Lambda^{2,1} \} \cap \Img \{ \del :
\Lambda^{1,0} \rightarrow \Lambda^{2,0} \} = \{0\}.
\end{align*}
\begin{proof} Let $\phi = \del
\ga$, $\delb \phi = 0$, $\ga \in \Lambda^{1,0}$.  A general calculation for
complex surfaces shows that for any metric $g$,
\begin{align*}
\brs{\phi}^2 dV_g =&\ \phi \wedge \bar{\phi}.
\end{align*}
Thus
\begin{align*}
\brs{\brs{\phi}}_{L^2}^2 = \int_M \phi \wedge \bar{\phi} = \int_M \phi \wedge
\delb \bar{\ga}= 0
\end{align*}
by Stokes theorem.
\end{proof}
\end{lemma}

\begin{lemma} \label{reductionlemma2} Let $(M^4, \omega, J)$ be a complex
surface
with pluriclosed
metric.  Suppose $b_1$ is odd and $h^{0,2} = 0$.  Then $[\del \omega] \neq 0 \in
H^3(M, \mathbb C)$ and $[\del \omega] \neq 0 \in H^{2,1}(M)$.
\begin{proof}  Suppose that $\del \omega$ is a $d$-exact form. 
By Lemma \ref{preliminarylemma} we conclude that
\begin{align*}
\del \omega = \delb \gb.
\end{align*}
Note this also holds trivially if we assume $[\del \omega] = 0 \in H^{2,1}(M)$.
Now, $\bar{\gb} \in \Lambda^{0,2}$ obviously satisfies $\delb \bar{\gb} = 0$. 
However, since $h^{0,2} = 0$ one can write
\begin{align*}
\bar{\gb} = \delb \bar{\ga}.
\end{align*}
Therefore $\del \omega = \delb \del \ga$ and, taking conjugates, $\delb \omega =
\del \delb \bar{\ga}$.  Let $\til{\omega} = \omega - \del \ga - \delb
\bar{\ga}$.  One computes directly that
\begin{align*}
d \til{\omega} =&\ \left( \del + \delb \right) \left( \omega - \del \ga - \delb
\bar{\ga} \right)\\
=&\ \del \omega - \delb \del \ga + \delb \omega - \del \delb \bar{\ga}\\
=&\ 0.
\end{align*}
Since the $(1,1)$ component of $\til{\omega}$ is positive definite, it follows
that
\begin{align*}
\int_M \til{\omega} \wedge \til{\omega} > 0.
\end{align*}
Since $b_1$ is odd, the intersection form of $M$ is negative definite
(\cite{BPV} Theorem IV.2.14), so this is a contradiction.
\end{proof}
\end{lemma}

\begin{thm} \label{reduction} Let $(M^{4}, \omega, J)$ be a complex surface
with pluriclosed metric
satisfying $h^{0,1} \leq 1$.  Let
\begin{align*}
B = \pi c_1(\omega) + \del \gb + \delb \bar{\gb} \in \pi c_1 \in \mathcal
H_{\del + \delb}.
\end{align*}
Suppose $\to = \omega + \del \ga + \delb \bar{\ga}$ is a solution to
\begin{align} \label{MAeqn}
\del \del^*_{\til{\gw}} \til{\gw} + \delb \delb^*_{\til{\gw}} \til{\gw} +
\frac{\sqrt{-1}}{2} \del \delb \log \det \til{g} = B.
\end{align}
Then
\begin{align} \label{reducedeqn}
\gb =&\ \del^*_{\to} \to + \frac{\sqrt{-1}}{4} \delb \log
\frac{\det{\til{g}}}{\det g}.
\end{align}

\begin{proof} Taking $\del$ of equation (\ref{MAeqn}) yields
\begin{align*}
0=&\ \del \delb \delb^*_{\to} \to - \del \delb \bar{\gb},
\end{align*}
or equivalently,
\begin{align*}
0 =&\ \delb \left[ \del \delb^*_{\to} \to - \del \bar{\gb} \right]
\end{align*}
Therefore $\del \delb^*_{\to} \to - \del \bar{\gb}$ is a $\delb$-closed,
$\del$-exact $(2,0)$-form.  Using Lemma \ref{reductionlemma} we conclude that
\begin{align*}
\del \delb^*_{\to} \to - \del \bar{\gb} = 0.
\end{align*}
Conjugating yields
\begin{align*}
\delb \left[\del^*_{\to} \to - \gb \right] = 0.
\end{align*}
Thus we may write the Hodge decomposition of $\del^*_{\to} \to - \gb$ with
respect to $\gD_{\delb, \gw}$ as
\begin{align*}
\del^*_{\to} \to - \gb =&\ h + \delb f
\end{align*}
where $h \in H^{0,1}$.  Next we claim that $h$ vanishes.  Once we know this,
differentiating and plugging into (\ref{MAeqn}) yields $f = -\frac{\sqrt{-1}}{4}
\delb \log \frac{\det \til{g}}{\det g}$ and the theorem follows.  First of all,
if $h^{0,1} = 0$ this is trivial, and since this observation holds in general
dimension we record this as Proposition \ref{generalreduction} below.  Next
suppose $h^{0,1} = 1$.  We compute
\begin{align*}
\int_M \left< h, \del^*_{\gw} \gw \right>_{\gw} dV =&\ \int_M \left<
\del^*_{\to} \to - \gb - \delb f, \del^*_{\gw} \gw \right>_{\gw} dV
\end{align*}
First observe
\begin{align*}
\int_M \left< \delb f, \del^*_{\gw} \gw \right>_{\gw} dV
=&\ \int_M \left< \del \delb f, \omega \right>_{\omega} dV_g\\
=&\ \int_M \del \delb f \wedge \omega\\
=&\ 0
\end{align*}
by Stokes Theorem, using that $\del \delb \omega = 0$.
Next we compute
\begin{align*}
\int_M \left< -\gb, \del^*_{\gw} \gw \right>_{\gw} dV =&\ \int_M \left< - \del
\gb, \gw \right>_{\gw} dV\\
=&\ \int_M \left< - \frac{1}{2} \left( \del \gb + \delb \bar{\gb} \right), \gw
\right>
\end{align*}
since $\gw$ is real.  Also we compute using (\ref{MAeqn}),
\begin{align*}
\int_M \left< \del^*_{\to} \to, \del^*_{\gw} \gw \right>_{\gw} dV =&\ \int_M
\left< \del \del^*_{\to} \to, \gw \right>_{\gw} dV\\
=&\ \int_M \left< \frac{1}{2} \left( \del \del^*_{\to} \to + \delb \delb^*_{\to}
\to \right), \gw \right>_{\gw} dV\\
=&\ \int_M \left< \frac{1}{2} \left( c_1(\to) - c_1(\omega) + \del \gb + \delb
\bar{\gb} \right), \gw \right>_{\gw} dV\\
=&\ \int_M \left< \frac{1}{2} \left( \del \gb + \delb \bar{\gb} \right), \gw
\right>_{\gw} dV
\end{align*}
Where in the last line we used that $c_1(\to) - c_1(\gw) = \del \delb \phi$ and
used that $\omega$ is orthogonal to the image of $\del\delb$ since it is
pluriclosed.  It follows that
\begin{align*}
\int_M \left< h, \del^*_{\omega} \omega \right>_{\gw} dV = 0
\end{align*}
However, by Lemma \ref{reductionlemma2},
$[\del \omega] \neq 0$, and $h^{2,1} = 1$ by Serre duality, therefore there is a
nonzero constant $a$ such that
$[\del \omega] = [a * h]$
(since $h$ is $\delb^*$-closed ,$* h$ defines a cohomology class).  Thus
\begin{align*}
0 =&\ \int_M \left< h, \del^*_{\gw} \gw \right>\\
=&\ \int_M \left< * h, \del \gw \right>\\
=&\ a \int_M \brs{h}^2
\end{align*}
therefore $h = 0$.
\end{proof}
\end{thm}

\begin{prop} \label{generalreduction} Let $(M^{2n}, \omega, J)$ be a complex
manifold with pluriclosed
metric satisfying $h^{0,1} = 0$, and
\begin{align*}
\Ker \{ \delb : \Lambda^{2,0} \rightarrow \Lambda^{2,1} \} \cap \Img \{ \del :
\Lambda^{1,0} \rightarrow \Lambda^{2,0} \} = \{0\}.
\end{align*}
Let
\begin{align*}
B = \pi c_1(\omega) + \del \gb + \delb \bar{\gb} \in \pi c_1 \in \mathcal
H_{\del + \delb}.
\end{align*}
Suppose $\to = \omega + \del \ga + \delb \bar{\ga}$ is a solution to
\begin{align}
\Phi(\to) = B.
\end{align}
Then
\begin{align*}
\gb =&\ \del^*_{\to} \to + \frac{\sqrt{-1}}{4} \delb \log
\frac{\det{\til{g}}}{\det g}.
\end{align*}
\begin{proof} The proof is clear from the proof of Theorem \ref{reduction}.
\end{proof}
\end{prop}

\noindent Theorem \ref{reduction} may be applied to reduce solutions to
(\ref{PCF}) on
surfaces.

\begin{thm} Let $(M^4, g_0, J)$ be a pluriclosed surface.  Say the solution to
(\ref{PCF}) with initial condition $g_0$ exists on a maximal time interval $[0,
T), T < \infty$.  Fix a background metric $\rho$ and express
\begin{align*}
\omega(t) = \omega_0 + \del \ga(t) + \delb \bar{\ga}(t) - t c_1(\rho).
\end{align*}
Then $\ga$ satisfies
\begin{align*}
\dt \ga =&\ \del^*_{\omega} \omega + \frac{\sqrt{-1}}{4} \delb \log \frac{\det
\til{g}}{\det g}
\end{align*}
\begin{proof} We differentiate the expression for $\omega(t)$ and use equation
(\ref{PCF}) to compute
\begin{align*}
\del \dot{\ga} + \delb \dot{\bar{\ga}} - c_1(\rho) = - \Phi(\omega).
\end{align*}
One may apply Theorem \ref{reduction} with $\gb = \dot{\ga}$ to conclude the
result.
\end{proof}
\end{thm}

\noindent Thus we have canonically reduced solutions to (\ref{PCF}) to solutions
of this
equation, coupled to an ODE.  There is a natural gauge to equation (\ref{aeqn}).

\begin{defn} Given $(M^{2n}, g, J)$ a complex manifold with Hermitian metric
$g$, let $\ga \in \Lambda^{0,1}$ and let $\psi \in \Lambda^{1,1}, d \psi = 0$. 
Let
\begin{gather}
[\ga, \psi] = \left\{ (\ga + \gb, \psi - \del \gb - \delb \bar{\gb}) | \gb \in
\Lambda^{0,1}, \delb \gb = 0 \right\}.
\end{gather}
We will refer to $[\ga, \psi]$ as the \emph{gauge equivalence class} of $(\ga,
\psi)$.
\end{defn}

\begin{prop} Let $(\ga(t), \psi(t))$ be a solution to (\ref{aeqn}).  Then
$(\ga(t), \psi(t))$ is gauge-equivalent to a pair $(\til{\ga}(t),
\til{\psi}(t))$ such that $\til{\ga}(t)$ solves a parabolic equation.
\begin{proof}  We find a gauge near time $t = 0$ for which $\til{\ga}(t)$ solves
a parabolic equation.  We take the $\del$-Hodge decomposition of $\ga(t)$ with
respect to $\omega_0$.  Specifically, consider
\begin{align*}
\ga(t) = \delb^*_{\omega_0} \nu(t) + \phi(t) + \delb f(t)
\end{align*}
where $\delb \phi(t) = \delb^*_{\omega_0} \phi = 0$.  Define
\begin{align*}
\til{\ga}(t) =&\ \ga(t) - \phi(0) - \delb f(0)\\
\til{\psi}(t) =&\ \psi(t) - \del \left( \phi(0) + \delb f(0) \right) - \delb
\left( \bar{\phi}(0) + \del \bar{f}(0) \right)
\end{align*}
Note that by construction $\delb^*_{\omega_0} \til{\ga}(t) = 0$.  Also, note
that $\til{\omega} = \omega_0 + \del \til{\ga} + \delb \bar{\til{\ga}} +
\til{\psi}$ by construction as well.  Furthermore, we have that
\begin{align*}
\dt \til{\ga}_{|t = 0} =&\ \del^*_{\omega_0} \omega_0 + \frac{\sqrt{-1}}{4}
\delb \log \frac{\omega_0^{ n}}{\rho^{n}}
\end{align*}
Our aim is to show that the right hand side is an elliptic operator for $\ga$. 
For a given metric $\omega$ we have the general coordinate formula
\begin{align*}
\left(\del^*_{\omega} \omega \right)_{\bj} =&\ \hi g^{p \bq} \left( \del_{\bq}
g_{p \bj} -
\del_{\bj} g_{p \bq} \right)
\end{align*}
Likewise we have
\begin{align*}
\left(\frac{\sqrt{-1}}{4} \delb \log \frac{{\omega}^{n}}{\rho^{
n}}\right)_{\bj}
=&\ \frac{\sqrt{-1}}{4} \left( g^{p \bq} \del_{\bj} g_{p \bq} - \rho^{p \bq}
\del_{j}
\rho_{p \bq} \right).
\end{align*}
Specializing these two formulas to the case $\til{\omega} = \omega + \del
\til{\ga} + \delb \bar{\til{\ga}} + \til{\psi}$ we compute
\begin{align*}
\left(\del^*_{\til{\omega}} \til{\omega} + \frac{\sqrt{-1}}{4} \delb \log
\frac{\til{\omega}^{n}}{\rho^{n}} \right)_{\bj} =&\
\frac{\sqrt{-1}}{2}
\til{g}^{p \bq} \left(  \del_{\bq} \left( \del_p \til{\ga}_{\bj} + \del_{\bj}
\til{\ga}_p
\right) - \del_{\bj} \left( \del_p \til{\ga}_{\bq} + \del_{\bq} \til{\ga}_p 
\right) \right)\\
&\ + \frac{\sqrt{-1}}{4} \til{g}^{p \bq} \left( \del_{\bj} \left( \del_p
\til{\ga}_{\bq} +
\del_{\bq} \til{\ga}_p \right) \right) + \mathcal O(\del \til{\ga})\\
=&\ \frac{\sqrt{-1}}{2} \til{g}^{p \bq} \del_p \del_{\bq} \til{\ga}_{\bj} +
\frac{\sqrt{-1}}{4} \til{g}^{p \bq} \del_{\bj} \del_{\bq} \til{\ga}_p -
\frac{\sqrt{-1}}{4} \til{g}^{p \bq} \del_{\bj} \del_{p} \til{\ga}_{\bq}.
\end{align*}
Now, using the condition $\delb^*_{\omega_0} \til{\ga} = 0$, we compute that
\begin{align*}
0 =&\ \del_{\bj} \delb^*_{\omega_0} \til{\ga}\\
=&\ g^{p \bq} \del_{\bj} \del_{p} \til{\ga}_{\bq} + \mathcal O(\del \til{\ga})
\end{align*}
Likewise since $\del^*_{\omega_0} \bar{\til{\ga}} = 0$ we have $g^{p \bq}
\del_{\bj} \del_{\bq} \til{\ga}_{p} = \mathcal O(\del \til{\ga})$.  Since
$\til{g}(0) = g_0$ we conclude that
\begin{align*}
\left(\del^*_{\til{\omega}} \til{\omega} + \frac{\sqrt{-1}}{4} \delb \log
\frac{\til{\omega}^{n}}{\rho^{n}} \right)_{\bj}(0) =&\
\frac{\sqrt{-1}}{2} \til{g}^{p \bq} \del_p \del_{\bq} \til{\ga}_{\bj}.
\end{align*}
which is a strictly elliptic operator.  The result follows.
\end{proof}
\end{prop}

\section{Pluriclosed flow as a gradient flow} \label{gradient}

In this section we exhibit that (\ref{PCF}) is the gradient flow of the first
eigenvalue of a certain Schr\"odinger operator associated to the time-dependent
metric.  What we actually show is that, after pulling back a solution to
(\ref{PCF}) by the one-parameter family of diffeomorphisms generated by the
vector field dual to the Lee form, one produces a solution to the
renormalization group flow of a nonlinear sigma model arising in string theory
(see \cite{Polchinski} 108-112).  This surprising fact both exhibits a
connection between pluriclosed flow and mathematical physics, and from another
point of view produces a large class of interesting examples of the
renormalization group flow.

Let us recall some notation from the introduction.  Fix $(M^{2n}, g, J)$ a
complex manifold with pluriclosed metric.  Let $\N$ denote the Bismut
connection, $\Rc$ the Ricci tensor of the Bismut connection, $\Rc^g$ the Ricci
curvature of $g$, $T$ the torsion of $\N$, and
\begin{align*}
\gt = - J d^* \omega
\end{align*}
the Lee form of $\omega$.  Lastly, as in the introduction, let $P$ denote the
representative of $c_1(M, J)$ associated to the Bismut connection $\N$.  We need
to show some identities relating these
tensors.  We start by recording some basic calculations which appear in
\cite{IvPap}.  It is important to remember below that $\Rc$ is not symmetric,
and $P$ is not $(1,1)$.

\begin{lemma} \label{lem1} (\cite{IvPap} Proposition 3.1) Let $(M^{2n}, g, J)$
be a pluriclosed structure.  Then
\begin{align*}
\Rc^g(X, Y) =&\ \Rc(X, Y) + \frac{1}{2} d^* T(X, Y) + \frac{1}{4} \sum_{i =
1}^{2n} g(T(X, e_i), T(Y, e_i))\\
P(X, Y) =&\ \Ric(X, JY) + \N_X \gt(J Y)\\
\Rc(Y, JX) + \Rc(X, JY) =&\ -(\N_X \gt)(JY) - \N_Y(\gt)(JX)\\
P(JX, JY) - P(X, Y) =&\ d^* T(JX, Y) - d^{\N} \theta(JX, Y)
\end{align*}
where $d^{\N}$ is the exterior derivative induced by $\N$.
\begin{proof} Note that the tensor $\gl$ from \cite{IvPap} vanishes when
$\del\delb \omega = 0$.  The third line is (\cite{IvPap} (3.9)).
\end{proof}
\end{lemma}

Let
\begin{align*}
H(X, Y) := P^{1,1}(JX, Y).
\end{align*}
In particular, note that (\ref{PCF}) is equivalent to
\begin{align} \label{metricflow}
\dt g =&\ - H.
\end{align}

\begin{lemma} \label{lem2} Let $(M^{2n}, g, J)$ be a pluriclosed structure. 
Then
\begin{align*}
H(X, Y) =&\ \frac{1}{2} \left[\Rc(X, Y) + \Rc(JX, JY) + \N_X \theta(Y) + \N_{JX}
\gt(JY) \right]
\end{align*}
\begin{proof} We directly compute using Lemma \ref{lem1}:
\begin{align*}
H(X, Y) =&\ P^{1,1}(JX, Y)\\
=&\ \frac{1}{2} \left[ P(JX, Y) + P(J J X, J Y) \right]\\
=&\ \frac{1}{2} \left[P(JX, Y) - P(X, JY)\right]\\
=&\ \frac{1}{2} \left[\Ric(JX, JY) + (\N_{JX} \gt)(JY) - \Ric(X, J J Y) - (\N_X
\gt)(J J Y)\right]\\
=&\ \frac{1}{2} \left[\Ric(X, Y) + \Ric(JX, JY) + (\N_X \gt)(Y) + (\N_{JX} \gt)
(JY)\right].
\end{align*}
\end{proof}
\end{lemma}

\begin{prop} \label{PCFgev} Given $(M^{2n}, \omega(t), J)$ a solution to
(\ref{PCF}), one has
\begin{align*}
\dt g =&\ \left[- \Rc^g + \frac{1}{4} \sum_{i = 1}^{2n} g(T(X, e_i), T(Y, e_i))
- \frac{1}{2} \mathcal L_{\gt^{\sharp}} g \right],
\end{align*}
where $\theta^{\sharp}$ is the vector field dual to $\theta$, taken with respect
to the time varying metric.
\begin{proof} Using the third line of Lemma \ref{lem1}, we compute
\begin{align*}
\Ric(JX, JY) + (\N_{JX} \gt)(JY) =&\ - \Ric(Y, JJX) - (\N_{Y} \gt)(JJX)\\
=&\ \Ric(Y, X) + (\N_Y \gt)(X).
\end{align*}
Plugging this into Lemma \ref{lem2} yields
\begin{align*}
2 H(X, Y) =&\ \Ric(X, Y) + \Ric(Y, X) + (\N_X \gt)(Y) + (\N_Y \gt)(X).
\end{align*}
The first two terms are twice the symmetric part of $\Ric$, which is easily
computed from the first line of Lemma \ref{lem1}.  It remains to show that the
last two terms are $\mathcal L_{\gt^{\sharp}} g$.  To do this we compute in
coordinates,
if $\Gamma$ denotes the connection coefficients of the Bismut connection and
$\Gamma^{LC}$ the Levi-Civita connection,
\begin{align*}
\N_i \gt_j =&\ \del_i \gt_j - \Gamma_{i j}^k \gt_k\\
=&\ \del_i \gt_j - \left( \Gamma^{LC} + \frac{1}{2} T \right)_{i j}^k \gt_k\\
=&\ D_i \gt_j - T_{i j}^k \gt_k
\end{align*}
where of course $D$ denotes the Levi-Civita derivative.  But $T$ is totally
skew, thus
\begin{align*}
\N_i \gt_j + \N_j \gt_i =&\ D_i \gt_j + D_j \gt_i - \frac{1}{2} T_{i j}^k \gt_k
- \frac{1}{2} T_{j i}^k \gt_k\\
=&\ D_i \gt_j + D_j \gt_i\\
=&\ \left(\mathcal L_{\gt^{\sharp}} g \right)_{ij}.
\end{align*}
\end{proof}
\end{prop}

\begin{prop} \label{PCFTev} Given $(M^{2n}, \omega(t), J)$ a solution to
(\ref{PCF}), one has
\begin{align*}
\dt T =&\ \frac{1}{2} \left[ \gD_{LB, g(t)} T - \mathcal L_{\gt^{\sharp}} T
\right].
\end{align*}
\begin{proof} Recall that $T = d^c \omega$, where $d^c = i(\delb - \del)$. 
Therefore
\begin{align*}
\dt T =&\ \dt d^c \omega\\
=&\ - d^c (P^{1,1}).
\end{align*}
Now note that, since $P$ is closed,
\begin{align*}
0 =&\ d P\\
=&\ (\del + \delb)(P^{1,1} + P^{2,0} + P^{0,2})\\
=&\ \del P^{1,1} + \delb P^{1,1} + \del P^{2,0} + \delb P^{2,0} +
\del P^{0,2} + \delb P^{0,2}. 
\end{align*}
By examining types we conclude from this the equations
\begin{align*}
0 =&\ \del P^{2,0} = \delb P^{0,2}\\
\del P^{1,1} =&\ - \delb P^{2,0}\\
\delb P^{1,1} =&\ - \del P^{0,2}.
\end{align*}
It follows that
\begin{align*}
- d^c (P^{1,1}) =&\ - i (\delb - \del) (P^{1,1})\\
=&\ i \del P^{0,2} - i \delb P^{2,0}.
\end{align*}
For convenience, set
\begin{align*}
\psi = d^* T - d^{\N} \theta.
\end{align*}
Now fix local complex coordinates, and compute using the last line of Lemma
\ref{lem1}
\begin{align*}
\left(i \del P^{0,2} \right)_{i \bj \bk} =&\ - \frac{i}{2} \del_i \left[ d^*
T(J \del_{\bj}, \del_{\bk}) - d^{\N} \theta(J \del_{\bj}, \del_{\bk}) \right]\\
=&\ - \frac{1}{2} \del_i \left[d^* T(\del_{\bj}, \del_{\bk}) - d^{\N}
\gt(\del_{\bj}, \del_{\bk})\right]\\
=&\ - \frac{1}{2} \left(\del \psi^{0,2} \right)_{i \bj \bk}.
\end{align*}
Likewise we can compute
\begin{align*}
\left(- i \delb P^{2,0} \right)_{\bi j k} =&\ \frac{i}{2} \del_{\bi} \left[
d^* T(J \del_{j}, \del_{k}) - d^{\N} \theta(J \del_{j}, \del_{k}) \right]\\
=&\ - \frac{1}{2} \del_{\bi} \left[d^* T(\del_{j}, \del_{k}) - d^{\N}
\gt(\del_{j}, \del_{k})\right]\\
=&\ - \frac{1}{2} \left(\delb \psi^{2,0} \right)_{\bi j k}.
\end{align*}

Note that it is a consequence of the last line of Lemma \ref{lem1} that
$\psi^{1,1} = 0$.  In particular, we have $\frac{1}{2} \psi = \rho^{2,0} +
\rho^{0,2}$ and so
\begin{align*}
\del \psi^{2,0} = \delb \psi^{0,2} = 0.
\end{align*}
Collecting these calculations yields
\begin{align*}
- d^c P^{1,1} =&\ -\frac{1}{2} d \psi\\
=&\ -\frac{1}{2} \left(d d^* T - d d^{\N} \theta \right).
\end{align*}
Since $T$ is closed, it follows that $d d^* T = - \gD_{LB, g(t)} T$.  Finally,
we observe a formula for $d^{\N} \theta$.
\begin{align*}
(d^{\N} \theta)_{ij} =&\ \N_i \theta_j - \N_j \theta_i\\
=&\ \del_i \theta_j - \left( \Gamma^{LC} + \frac{1}{2} T \right)_{i j}^k
\theta_k - \del_j \theta_i + \left( \Gamma^{LC} + \frac{1}{2} T \right)_{j i}^k
\theta_k\\
=&\ \left(d \theta - \theta^{\sharp} \hook T \right)_{i j}
\end{align*}
It follows from the Cartan formula and the fact that $T$ is closed that
\begin{align*}
d d^{\N} \theta =&\ d \left( d \theta - \theta^{\sharp} \hook T \right)\\
=&\ - d \left( \theta^{\sharp} \hook T \right)\\
=&\ - \mathcal L_{\theta^{\sharp}} T + \theta^{\sharp} \hook \left(d T \right)\\
=&\ - \mathcal L_{\theta^{\sharp}} T.
\end{align*}
Therefore
\begin{align*}
- d^c P^{1,1} = \frac{1}{2} \left[\gD_{LB, g(t)} - \mathcal
L_{\theta^{\sharp}} T \right].
\end{align*}
\end{proof}
\end{prop}

\begin{thm} \label{gaugefix} Let $(M^{2n}, \til{\omega}(t), J)$ be a solution to
(\ref{PCF}).  Let $X(t) = \frac{1}{2} \til{\theta}^{\sharp}$, where $\sharp$
means the vector dual taken with respect to the time-varying metric, and let
$\phi_t$ denote the one parameter family of diffeomorphisms generated by $X(t)$.
 Let $\til{T}$ denote the torsion of the time-varying Bismut connections.  Let
$(g(t), T(t)) = (\phi^*(\til{g})(t)), \phi_t^*(\til{T})(t))$.  Then
\begin{gather} \label{CRF}
\begin{split}
\dt g =&\ - \Rc^g + \frac{1}{4} \mathcal H\\
\dt T =&\ \frac{1}{2} \gD_{LB} T.
\end{split}
\end{gather}
where $\mathcal H_{ij} = g^{k l} g^{mn} T_{i k m} T_{j l n}$.
\begin{proof} This follows from a standard calculation using Propositions
\ref{PCFgev} and \ref{PCFTev}.
\end{proof}
\end{thm}

As noted above, the system of equations (\ref{CRF}) arises naturally in physics
as the renormalization group flow of a nonlinear sigma model.  By extending
Perelman's energy functional (\cite{Perelman}) to this coupled system, Oliynyk,
Suneeta, and Woolgar showed that (\ref{CRF}) is the gradient flow of a nonlinear
Schr\"odinger operator (\cite{Woolgar}).  To discuss this let us generalize the
notation slightly.  As in the introduction, let $(M^n, g)$ be a Riemannian
manifold, and let $T$ denote
a three-form on $M$.  Let
\begin{align*}
\mathcal F(g, T, f) =&\ \int_M \left[ R - \frac{1}{12} \brs{T}^2 + \brs{\N f}^2
\right] e^{-f} dV.
\end{align*}
Furthermore set
\begin{align*}
\gl(g, T) = \inf_{\{f | \int_M e^{-f} dV = 1 \}} \mathcal F(g, T, f).
\end{align*}

\begin{prop} \label{woolgar} (\cite{Woolgar} Proposition 3.1)  The gradient flow
of $\gl$ is
\begin{gather} \label{CRF2}
\begin{split}
\dt g =&\ - 2 \Rc + \frac{1}{2} \mathcal H - 2 \N^2 f,\\
\dt T =&\ \gD_{LB} T - d (\N f \hook T),
\end{split}
\end{gather}
where $f$ satisfies the conjugate heat equation
\begin{align} \label{conjheat}
\dt f =&\ - \gD f - R + \frac{1}{4} \brs{T}^2.
\end{align}
\end{prop}

\noindent For concreteness sake we now record the proof of Theorem
\ref{gradientPCF}.

\begin{proof}  Clearly equation (\ref{CRF2}) is diffeomorphism equivalent to
(\ref{CRF}).  By combining Proposition \ref{woolgar} with Theorem
\ref{gaugefix}, we obtain the statement of Theorem \ref{gradientPCF}.
\end{proof}

Furthermore, in \cite{FIL} Feldman, Ilmanen and Ni gave a generalization of
Perelman's steady and shrinking entropies to an entropy modeled on expanding
solitons.  Surprisingly, this expanding entropy has an extension to (\ref{CRF}),
as shown by the first named author.  Define

\begin{align*}
\mathcal W_+(g, T, u, \gs) =&\ \int_M \left[ \gs \left(\frac{\brs{\N u}^2}{u} +
R u - \frac{1}{12} \brs{T}^2 u \right) + u \log u \right] dV\\
=&\ \int_M \left[ \gs \left( \brs{\N f_+}^2 + R - \frac{1}{12} \brs{T}^2 \right)
- f_+ + n \right] u dV
\end{align*}
where $f_+$ is implicitly defined by
\begin{align*}
u = \frac{e^{-f_+}}{(4 \pi \gs)^{\frac{n}{2}}}.
\end{align*}

\begin{thm} \label{expandingentropy} (\cite{Streets} Theorem 6.2) Let $(M^n,
g(t), T(t))$ be a solution to (\ref{CRF}) on $[t_1, t_2]$ and suppose
$u(t)$ is the solution to (\ref{conjheat}).  Let
\begin{align*}
v_+ = \left[(t - t_1)(2 \gD f_+ - \brs{\N f_+}^2 + R - \frac{1}{12}
\brs{T}^2) - f_+ + n \right] u.
\end{align*}
Then
\begin{align*}
& \left( \dt + \gD - R + \frac{1}{4} \brs{T}^2 \right) v_+\\
&\ \qquad = 2 (t - t_1) \left( \brs{\Rc - \frac{1}{4} \mathcal H + \N^2 f_+ +
\frac{g}{2 t}}^2 + \frac{1}{4} \brs{d^* T - \N f_+ \hook T}^2 \right) u +
\frac{1}{6} \brs{T}^2 u.
\end{align*}
\end{thm}

\begin{cor} Let $(M^n, g(t), T(t))$ be a solution to (\ref{CRF}) on $[t_1,
t_2]$ and suppose $u(t)$ is a solution to the conjugate heat equation.  Then
\begin{align*}
\dt \mathcal W_+(g(t), T(t), u(t), t - \tau_1) =&\ \int_M 2 u \left[ (t -
t_1) \brs{\Rc - \frac{1}{4} \mathcal H + \N^2 f_+ + \frac{g}{2(t - t_1)}}^2
\right.\\
&\ \left. + \frac{1}{4} (t - t_1) \brs{ d^* T - \N f_+ \hook T}^2 +
\frac{1}{12} \brs{T}^2
\right] dV.
\end{align*}
\end{cor}

We can derive further corollaries from these results, akin to the ``ruling out
of breathers'' statements discovered by Perelman (\cite{Perelman}).  First
recall two definitions.

\begin{defn} We say that a solution to (\ref{CRF}) is a \emph{breather} if there
are times $t_1 < t_2$, a constant $\ga > 0$ and a diffeomorphism $\phi$ such
that $\ga g(t_1) = \phi^* g(t_2)$.  The breather is \emph{steady},
\emph{shrinking} or \emph{expanding} if $\ga = 1, \ga < 1$, or $\ga > 1$,
respectively.
\end{defn}

\begin{defn} We say that a solution to (\ref{CRF}) is a \emph{gradient soliton}
if there is a function $f$ and a constant $\gl$ so that 
\begin{align*}
0=&\ \Rc - \frac{1}{4} \mathcal H + \N^2 f - \gl g\\
0 =&\ \gD_{LB} T - d (\N f \hook T)
\end{align*}
The soliton is \emph{steady}, \emph{shrinking} or \emph{expanding} if $\gl = 0,
\gl > 0$, or $\gl < 0$, respectively.
\end{defn}

\begin{cor} \label{nobreathers} Any solution to (\ref{CRF}) which is a steady
breather is a steady soliton.  Any solution to (\ref{CRF}) which is an expanding
breather is an Einstein metric with $T \equiv 0$.
\begin{proof} The first statement follows immediately from Proposition
\ref{woolgar}.  For the second, we note that Theorem \ref{expandingentropy}
clearly implies that an expanding breather is an expanding soliton, and moreover
$T \equiv 0$.  Thus $g(t)$ is an expanding \emph{Ricci} soliton, which are known
to be negative constant Einstein metrics, a result originally due to Hamilton
(\cite{Hamilton2}).
\end{proof}
\end{cor}

\section{Nonsingular solutions} \label{nonsingularsection}

In this section we derive a strong topological consequence of the conjectural
regularity
picture of solutions to (\ref{PCF}) by ruling out nonsingular solutions of
(\ref{PCF}) on Class
$\seven^+$ surfaces.

\begin{thm} Suppose Conjecture \ref{strongexistenceconjn2} holds.  Then any
Class $\seven^+$ surface contains an irreducible effective divisor of
nonpositive
self intersection.
\begin{proof} We want to examine the volume-normalized version of (\ref{PCF}). 
Let
\begin{align*}
\psi(\omega) = \frac{\int_M \tr_{\omega} \left[ \del \del^*_{\gw} \gw + \delb
\delb^*_{\gw} \gw + \frac{\sqrt{-1}}{2} \del\delb \log \det g \right] dV}{\int_M
dV}.
\end{align*}
The volume normalized pluriclosed flow is
\begin{gather} \label{vflow}
\begin{split}
\dt \omega =&\ \del \del^*_{\gw} \gw + \delb \delb^*_{\gw} \gw +
\frac{\sqrt{-1}}{2} \del\delb \log \det g - \frac{1}{n} \psi \omega\\
\omega(0) =&\ \omega.
\end{split}
\end{gather}

Let $(M^4, J)$ be a complex surface of Class $\seven^+$.  Note that by the
theorem
of Gauduchon \cite{Gauduchon}, there are always pluriclosed metrics on complex
surfaces, so we can find an initial condition for (\ref{vflow}).  If $M$
contains no irreducible divisor of nonpositive self intersection, then
Conjecture \ref{strongexistenceconjn2} automatically implies that the solution
to (\ref{vflow}) with any initial condition exists for all time with a uniform
bound on curvature which moreover persists after scaling the diameter to unit
size.  Note that this implies the diameter is in fact bounded, for if not by
rescaling the diameter to unit size we would produce a sequence of metrics with
bounded curvature and volume approaching zero, which by \cite{CFG}, \cite{CGT}
would force
$\chi(M) = 0$, but $\chi(M) = b_2 > 0$.

We want to derive a contradiction from the existence of such a flow.  To do this
we first identify the qualitative behavior of the
corresponding solution to (\ref{PCF}).  Specifically, using monotone
quantities we see that this solution exists for all time with volume growing
quadratically.  Thus to obtain the solution to (\ref{vflow}) we must be
uniformly scaling
down this metric, and we will finish the proof by applying the expanding entropy
formula.  We begin with a
definition and a series of lemmas. 

\begin{defn} Let $(M^{2n}, g, J)$ be a Hermitian manifold.  Let the degree of
$(M, g)$ be
\begin{gather} \label{ddef}
d = \deg(M, g) := \int_M \left< c_1(M), \omega \right> = \int_M \left( -
\frac{\sqrt{-1}}{2}
\del
\delb \log \det g \right) \wedge \omega^{n-1}.
\end{gather}
More generally, given $\LL$ a line bundle over $M$, define
\begin{gather} \label{degdef}
\deg(\LL) := \int_M c_1(\LL) \wedge \omega^{n-1}.
\end{gather}
\end{defn}

\noindent Note that the definition of degree is typically made with respect to a
fixed Gauduchon metric, i.e. a metric satisfying $\del \delb \omega^{n-1} = 0$,
so that the value does not depend on the representative of $c_1$.  In the case
$n = 2$ Gauduchon metrics are the same as pluriclosed metrics, and the evolution
of the degrees of line bundles is particularly clean.

\begin{lemma} \label{degreeev} Let $(M^4, g(t), J)$ be a solution to (\ref{PCF})
on a complex
surface, and let $L$ be a line bundle over $M$.  Then
\begin{align*}
\dt \deg_{g_t}(L) =&\ - c_1(L)\cdot c_1(M).
\end{align*}
\end{lemma}

\begin{lemma} \label{volumeev} Let $(M^{4}, g(t), J)$ be a solution to
(\ref{PCF}).  Then the volume of $g(t)$ satisfies
\begin{align*}
\dt \Vol(g(t)) =&\ 2 \int_M \brs{\del^* \omega}^2 - d.
\end{align*}
\end{lemma}

\noindent Next we would like to specialize these to the case of Class $\seven$
surfaces.

\begin{lemma} \label{class7degreeev} Let $(M^4, g(t), J)$ be a solution to
(\ref{PCF}) on a Class $\seven$ surface with $b_2 = n$.  Then
\begin{align*}
\deg_{g(t)}(M) =&\ \deg_{g(0)}(M) + n t.
\end{align*}
\begin{proof} This follows immediately from Proposition \ref{degreeev} and the
fact that for Class $\seven$ surfaces, $c_1^2 = -n$.
\end{proof}
\end{lemma}

\begin{prop} \label{unnormprop1} Let $(M^4, J)$ be a compact Class $\seven^+$
surface.  Suppose
${\gw}(t)$ is a solution to (\ref{vflow}) on $M \times [0,\infty)$ with
uniformly bounded curvature.  Then the
corresponding solution
to (\ref{PCF}) exists on $[0, \infty)$.
\begin{proof}  Suppose the corresponding solution to unnormalized flow existed
on $[0, \tau), \tau < \infty$.  First note that the degree of $M$ remains finite
on
$[0, \tau)$.  However, to rescale to get the volume normalized flow we must be
rescaling by a factor going to infinity since the curvature must be blowing up,
and we have assumed the volume-normalized flow is nonsingular.  Thus the volume
must be going to zero.  Using Lemma \ref{volumeev}
we see that at some point the degree must be positive.  But this condition is
preserved, since the degree grows linearly by Lemma \ref{class7degreeev}.  In
taking the rescaling limit,
this says that the degree must go to infinity as $t \rightarrow \infty$ in the
volume
normalized flow.  Since by assumption the volume normalized equation has bounded
curvature this is a contradiction.
\end{proof}
\end{prop}

\begin{prop} \label{unnormprop2} Let $(M^4, J)$ be a compact Class $\seven^+$
surface.  Suppose
$\til{\omega}(t)$ is a solution to (\ref{vflow}) on $M \times [0,\infty)$ with
uniformly bounded curvature.  Then if
$\omega(t)$ is the corresponding solution
to (\ref{PCF}), there exists a constant $C$ such that
\begin{align*}
\frac{1}{C}\left(1 + t^2 \right) \leq \Vol(g(t)) \leq C \left(t^2 + 1 \right).
\end{align*}
\begin{proof}  By assumption the corresponding solution to (\ref{vflow}) has
bounded curvature, and of course bounded volume.  It follows that the
scale-invariant quantity
$\frac{d}{\Vol^{\frac{1}{2}}}$ is bounded along the solution to (\ref{vflow}). 
Thus this quantity is bounded along (\ref{PCF}) as well.  By Lemma
\ref{class7degreeev} it follows that
\begin{align*}
d(0) + n t =&\ d(t) \leq C \Vol(g(t))^{\frac{1}{2}}
\end{align*}
with $n > 0$.  The lower bound for $\Vol(g(t))$ follows by squaring the above
inequality.  Next
we note the evolution equation for the degree under (\ref{vflow}).  In
particular one has
\begin{align*}
\dt \deg(M) =&\ - c_1^2 + \deg(M)^2 - 2 \deg M \int_M \brs{\del^* \omega}^2.
\end{align*}
Since $\int_M \brs{\del^* \omega}^2$ is bounded, it easy follows that $\lim_{t
\rightarrow \infty} \deg(M) > \ge > 0$ for some $\ge > 0$.  It follows that the
scale invariant quantity
$\frac{\deg(M)}{V^{\frac{1}{2}}} \geq \ge > 0$.  Thus this inequality holds for
the unnormalized flow as well, hence
\begin{align*}
d(0) + n t =&\ d(t) \geq \ge \Vol(g(t))^{\frac{1}{2}}.
\end{align*}
The upper volume bound now follows, completing the proof.
\end{proof}
\end{prop}

\noindent We now give the proof of the theorem.  From the proposition above we
see that the solution to (\ref{vflow})
is uniformly equivalent to a solution of 
\begin{gather} \label{vflow2}
\begin{split}
\dt \omega =&\ \del \del^*_{\gw} \gw + \delb \delb^*_{\gw} \gw +
\frac{\sqrt{-1}}{2} \del\delb \log \det g - \omega\\
\omega(0) =&\ \omega.
\end{split}
\end{gather}
Moreover, this solution has uniformly bounded curvature and diameter.  We claim
that there
exists a uniform lower bound on the injectivity radius as well.  If there exists
a sequence of points $(x_i, t_i)$ with $t_i \rightarrow \infty$ such that
$\inj_{g_{t_i}}(x_i) \rightarrow 0$, then since there is a uniform diameter
bound it follows from \cite{CGT} that $\inj_{g_{t_i}}(x) \rightarrow 0$ for all
points $x \in M$.  In particular, the manifold $M$ admits a sequence of metrics
collapsing with bounded curvature, which by \cite{CFG} implies that $\chi(M)
= 0$.  But for $(M^{4}, J)$ a complex Class $\seven^+$ surface one has $\chi(M)
= b_2(M) > 0$, so this is a contradiction, and so the lower injectivity radius
bound follows.

Thus we can construct a blowdown limit for the unnormalized flow.  Let $\gl_i
\rightarrow
\infty$ and set
\begin{align*}
g_i(t) := \frac{1}{\gl_i} g( \gl_i t),
\end{align*}
defined for $t > \frac{1}{\gl_i}$.  By the estimates we have shown and the
compactness theorem of (\cite{Streets}) there is a subsequence of $\{(M, g_i(t),
J) \}$, also denoted by which in the Cheeger-Gromov sense to a limiting
pluriclosed flow $\{(M_{\infty}, g_{\infty}(t), J_{\infty})\}$.  Since $\mathcal
W_+$ is
invariant under the blowdown rescaling, and is monotone increasing and bounded
above, it is clear that $W_+$ is constant along $g_{\infty}$, which thus must
be an expanding soliton, and hence is a K\"ahler-Einstein metric.  In
particular, we have that $(M_{\infty}, g_{\infty}(1), J)$ is a K\"ahler
manifold.  By the Hodge decomposition since $(M_{\infty}, J_{\infty})$ is
K\"ahler we have that $b_1(M_{\infty})$ is even.  On the other hand by the
diameter and injectivity radius bounds we have that $M \cong M_{\infty}$, and
$b_1(M) = 1$.  This is a contradiction, finishing the proof.
\end{proof}
\end{thm}

By general theory (\cite{Nakamura} Lemma 2.2) the curve is either a rational
curve, rational curve with double point, or an elliptic curve.  If the curve is
elliptic, the manifold is known (Nakamura \cite{Nakamura}, Enoki \cite{Enoki}). 
Furthermore, Class $\seven^+$
surfaces which contain $b_2$ rational curves automatically contain a global
spherical shell by the result of Dloussky, Oeljeklaus and Toma \cite{DOT}. 
Therefore we see that Conjecture \ref{strongexistenceconjn2} implies the
classification of Class $\seven^+$ surfaces with $b_2 = 1$, a theorem obtained
by
Teleman using gauge theory \cite{Tel2}.  Furthermore it implies a concrete
complex analytic conclusion on \emph{any} Class $\seven^+$ surface.  It seems
likely that a more detailed analysis of the limit points can yield the entire
classification of Class $\seven^+$ surfaces as a consequence of Conjecture
\ref{strongexistenceconjn2}.

\section{Conclusion} \label{conclusion}

Given the results contained herein, equation (\ref{PCF}) clearly seems to be a
very natural parabolic equation on complex manifolds.  By the results of
\cite{ST2}, the corresponding elliptic (static) equation, i.e.
\begin{align} \label{static}
P^{1,1} = \gl \omega
\end{align}
 seems very closely
related to the K\"ahler-Einstein condition, and we further seen here the
relationship of solutions to (\ref{PCF}) and the topological and complex
structure of surfaces.  While there are only a few large classes of examples of
complex manifolds of dimension $n \geq 3$ admitting pluriclosed (but not
K\"ahler) geometries, it seems likely that understanding the existence problem
for static metrics in higher dimensions will have relevance.  We take the time
here to observe some further structural results for static metrics in any
dimension.  First of all, we recall the Bochner formula for holomorphic forms on
complex manifolds.

\begin{thm} \label{bochner} (\cite{KWu}, \cite{Bochner}) Let $(M^{2n}, g, J)$ be
a Hermitian manifold.  Fix $\eta$ a holomorphic $(p,0)$-form.  Then
\begin{align} \label{boch}
\gD \brs{\eta}^2 =&\ \brs{\N \eta}^2 + \brs{\bar{\N} \eta}^2 + \left< S \circ
\eta, \eta \right>
\end{align}
here $\gD = \tr_{\omega} \del \delb$ is the canonical Laplacian and $\N$ is the
Chern connection.  Also, $S \circ \eta$ is the natural action induced on
$\Lambda^{p,0}$ of and endomorphism of $T^{1,0}$.  In particular, in
coordinates,
\begin{align*}
\left( S \circ \eta \right)_{i_1 \dots i_p} =&\ \frac{1}{p!} \sum_{j = 1}^p
S_{i_j}^k \eta_{i_1 \dots i_{j-1} k i_{j+1} \dots i_p}.
\end{align*}
\end{thm}

\begin{cor} (\cite{IvPap} Corollary 4.4) 
\begin{itemize}
\item {Let $(M^{2n}, \omega, J)$ be a compact
complex manifold with $\pi c_1 = 0$.  Suppose $\omega$ is a
static metric, which
necessarily has $s \equiv 0$.  Then every holomorphic $(p,0)$-form is parallel
with respect to the Chern connection.}
\item{ Let $(M^{2n}, \omega, J)$ be a compact complex manifold with $\pi c_1 >
0$.  Suppose $\omega$ is a static metric, which necessarily has $s \equiv c >
0$.  Then $H^0(M, \Lambda^p) = 0$, $p = 1, \dots n$.}
\end{itemize}
\end{cor}
Note that in the second part of the above corollary, we mean $c_1 > 0$ as a
class in the Aeppli cohomology group $\mathcal H^{1,1}_{\del + \delb}$ as
defined in section \ref{conesection}.  These are precisely the corollaries to
his conjecture observed by Calabi
\cite{Calabi} in the
K\"ahler setting, indicating that (\ref{PCF}) is a very natural extension of the
Calabi-Yau/K\"ahler Ricci flow theory.  We close with a final vanishing result
for static metrics

\begin{prop} \label{vanishinglemma} Let $(M^{2n}, g, J)$ be a complex manifold
with either $c_1 = 0$ or $c_1 > 0$, with $g$ a static metric.  Then either $g$
is K\"ahler or $h^{n-1,0} = 0$.
\begin{proof} Since $Q^1 \geq 0$, and $\gl \geq 0$ by the assumption on the
first Chern class, equation (\ref{static}) clearly implies that $S \geq Q^1$. 
Therefore by applying the maximum principle to (\ref{boch}) we conclude that
every holomorphic section $\eta$ of $\Lambda^{n-1,0}$ is parallel with respect
to the Chern connection.  In particular, it is of constant norm.  If $g$ is not
K\"ahler, there is a point $p \in M$ where the torsion tensor does not vanish
identically.  Specifically, we can pick complex coordinates where $S$, and hence
$Q^1$, are diagonalized.  Without loss of generality $T_{1 2 \bar{j}} \neq 0$. 
Thus $Q^1_{1 \bar{1}} > 0$, $Q^1_{2 \bar{2}} > 0$ (see Lemma
\ref{Phicoordinates} for the expression of $Q^1$).   It follows that $S(p) \geq
Q^1(p)$ is $n-1$ positive at $p$.  By the form of the Bochner formula
(\ref{boch}), if $\eta$ does not vanish we conclude that $\gD \brs{\eta}^2(p) >
0$, contradicting that $\eta$ is parallel.
\end{proof}
\end{prop}

\bibliographystyle{hamsplain}

\end{document}